\documentclass[12pt]{amsart} 
\usepackage{fancyhdr}

\usepackage{amsmath,amsthm,amscd,amsfonts,amssymb,euscript}

\usepackage[all]{xypic}
\usepackage{latexsym}
\usepackage{color} 
\theoremstyle{plain}

\begin{document}
\input{amssym.def}

\numberwithin{equation}{section}

\newtheorem{guess}{\sc Theorem}[section]
\newcommand{\bth}{\begin{guess}$\!\!\!${\bf }~}
\newcommand{\eeth}{\end{guess}}

\newtheorem{propo}[guess]{\sc Proposition}

\newcommand{\bprop}{\begin{propo}$\!\!\!${\bf }~}
\newcommand{\eprop}{\end{propo}}

\newtheorem{lema}[guess]{\sc Lemma}
\newcommand{\blem}{\begin{lema}$\!\!\!${\bf }~}
\newcommand{\elem}{\end{lema}}

\newtheorem{defe}[guess]{\sc Definition}
\newcommand{\bdefe}{\begin{defe}$\!\!\!${\bf }~}
\newcommand{\edefe}{\end{defe}}

\newtheorem{coro}[guess]{\sc Corollary}
\newcommand{\bcor}{\begin{coro}$\!\!\!${\bf }~}
\newcommand{\ecor}{\end{coro}}

\newtheorem{rema}[guess]{\it Remark}
\newcommand{\brem}{\begin{rema}$\!\!\!${\it }~\rm}
\newcommand{\erem}{\end{rema}}

\theoremstyle{remark}
\newtheorem{note}[guess]{Notation}
\newtheorem{example}{Example}[section]
\newcommand{\beg}{\begin{example}$\!\!\!${\bf }~~\rm}
\newcommand{\eeg}{\end{example}}
\theoremstyle{remark}
\newtheorem{observe}[guess]{Observation}
\newtheorem{assump}[guess]{Assumption}
\newtheorem{sub}[guess]{}{\bf }

\newtheorem{ack}{Acknowledgments}       
\renewcommand{\theack}{} 

\newcommand{\pf}{{\noindent \sc Proof. }}
\newcommand{\enpf}{\begin{flushright} {\it Q.E.D} \end{flushright}}

\newtheorem{eq}[guess]{equation}
\newcommand{\beq}{\begin{equation}}
\newcommand{\eeq}{\end{equation}}

\newtheorem{eqa}[guess]{eqnarray}
\newcommand{\beqa}{\begin{eqnarray}}
\newcommand{\eeqa}{\end{eqnarray}}

\newcommand{\G}{\Gamma}
\newcommand{\hra}{\hookrightarrow}
\newcommand{\lr}{\longrightarrow}
\newcommand{\llr}{\longleftarrow}
\newcommand{\K}{\mathrm K}
\newcommand{\mk}{\mathrm k}
\newcommand{\Ha}{\mathcal H}
\newcommand{\Rss}{\mathcal R^\Gamma}
\newcommand{\Gt}{\tilde \mathcal F ^ \Gamma}
\newcommand{\U}{\mathcal U}
\newcommand{\N}{\mathcal N}
\newcommand{\Ga}{\mathcal G}
\newcommand{\m}{\mathcal}
\newcommand{\g}{\goth }

\newcommand{\LL}{L^{(\alpha_1)}}
\newcommand{\Ld}{L^{(\alpha_1)*}}
\newcommand{\M}{M^{(\alpha_2)}}
\newcommand{\Md}{M^{(\alpha_2)*}}

\newcommand{\Vls}{\mathcal V \mid_{q_{\em l}^{-1}(\underline s)}}
\newcommand{\Vms}{\mathcal V \mid_{q_{\em m}^{-1}(\underline s)}}
\newcommand{\Vlas}{\mathcal V \mid_{q_{(\underline l)}^{-1}(\underline s)}}
\newcommand{\Vmas}{\mathcal V \mid_{q_{((\underline m))}^{-1}(\underline s)}}
\newcommand{\Vdms}{\mathcal V \mid_{q_{\em m}^{-1}(\underline s)}^{**}}
\newcommand{\Vym}{\mathcal V \mid_{C_{\em m}}}
\newcommand{\Vyl}{\mathcal V \mid_{C_{\em l}}}

\newcommand{\Sma}{S_{m}^\Gamma}
\newcommand{\Sm}{S_{\em m}^\Gamma}
\newcommand{\Sl}{S_{\em l}^\Gamma}
\newcommand{\Sal}{S_{(\underline l)}^\Gamma}
\newcommand{\Sam}{S_{(\underline m)}^\Gamma}
\newcommand{\Smi}{S_{\em m_i}^{\Gamma}}

\newcommand{\Cm}{C_{\em m}^{\Gamma}}
\newcommand{\Cn}{C_{\em {m_0}}^{\Gamma}}
\newcommand{\Cam}{C_{\em (m)}^{\Gamma}}

\newcommand{\Cal}{C_{\em (l)}^{\Gamma}}

\newcommand{\Km}{K_{\em m}^\Gamma}

\newcommand{\un}{\em}
\newcommand{\us}{\em s}

\newcommand{\qu}{q_{\em m}}
\newcommand{\qm}{q_{{\em m}_0}}

\newcommand{\pu}{p_{\em m}}
\newcommand{\phim}{\phi_{\em m}}

\newcommand{\gb}{$\Gamma$-bundle}
\newcommand{\gsb}{$\Gamma$-subbundle}

\newcommand{\qis}{q_{\em m}^{-1}(\underline s)}

\newcommand{\Ul}{U_{\em l}}
\newcommand{\Um}{U_{\em m}}
\newcommand{\Uam}{U_{(\underline m)}}
\newcommand{\Ual}{U_{(\underline l)}}

\newcommand{\Wm}{\mathcal W_{\em m}}

\newcommand{\Gm}{\tilde{\mathcal G_m}}
\newcommand{\Gn}{\tilde{\mathcal G_{m_0}}}
\newcommand{\pmV}{\phi_{\em m}^*(p_{\em m}^* V)}

\newcommand{\wpmV}{\mathbb P (\wedge^{r}\phi_{\em
    m}^*(p_{\em m}^* V))}

\newcommand{\pmwV}{\mathbb P (\phi_{\em m}^*(p_{\em
    m}^*(\wedge^{r}V)))}

\newcommand{\pV}{\phi_{\em {m_0}}^*(p_{\em {m_0}}^* V)}

\newcommand{\E}{\mathcal E}
\newcommand{\Ep}{\mathcal E'}
\newcommand{\Epp}{\mathcal E''}

\newcommand{\bq}{{\mathbb Q}}
\newcommand{\cc}{{\mathcal C}}
\newcommand{\cu}{{\mathcal U}}
\newcommand{\cf}{{\mathcal F}}
\newcommand{\ce}{{\mathcal E}}
\newcommand{\co}{{\mathcal O}}
\newcommand{\cg}{{\mathcal G}}
\newcommand{\cm}{{\mathcal H}}
\newcommand{\cn}{{\mathcal N}}

\newcommand{\mcp}{{\mathcal P}}
\newcommand{\bp}{{\mathbb P}}
\newcommand{\cp}{{\sf P}}
\newcommand{\bpp}{{\it ps}}

\newcommand{\cq}{{\sf Q}}
\newcommand{\cs}{{\mathcal S}}
\newcommand{\cl}{{\mathcal L}}

\newcommand{\ctext}[1]{\makebox(0,0){#1}}
\setlength{\unitlength}{0.1mm}

\def\bC{{\Bbb C}}
\def\a{\alpha }
\def\bR{{\Bbb R}}
\def\tX{{\tilde X}}

\newcommand{\Mo}{M_\Gamma^{\mu ss}}
\newcommand{\Ro}{\R_\Gamma^{\mu ss}}
\newcommand{\Q}{Quot(\mathcal H, P)}
\newcommand{\gF}{gr^\mu_\G(F)}
\newcommand{\Fs}{gr^\mu_\G(F)^*}

\title[Stable bundles]{An analogue of the Narasimhan-Seshadri theorem in higher dimensions and some applications}
\author{V. Balaji}{\thanks{Research of the first author was partly supported by the J.C. Bose Research grant.}} 
\address{Chennai Math.Inst.SIPCOT IT Park, Siruseri-603103, India,
balaji@cmi.ac.in}
\author{A.J. Parameswaran}
\address{Kerala School of Mathematics, Kozhikode, Kerala and School of Mathematics, Tata Inst. of Fund. Research, Mumbai-400095, India, param@math.tifr.res.in}
\date{}
\dedicatory {To Madhav Nori on his sixtieth birthday}

\begin{abstract}
We prove an analogue in higher dimensions of the classical Narasimhan-Seshadri theorem for strongly stable vector bundles of degree $0$ on a smooth projective variety $X$ with a fixed ample line bundle $\Theta$. As applications, over fields of characteristic zero, we give a new proof of the main theorem in a recent paper of Balaji and Koll\'ar  and derive an effective version of this theorem;  over uncountable fields of positive characteristics, if $G$ is a simple and simply connected algebraic group and the characteristic of the field is bigger than the Coxeter index of $G$, we prove the existence of 
strongly stable principal $G$ bundles on smooth projective surfaces whose holonomy group is the whole of $G$.
\end{abstract}

\maketitle
\tiny
\tableofcontents
\normalsize
\section{Introduction}
Let $X$ be a smooth projective variety over an algebraically closed field $k$ of arbitrary characteristic. Let $\Theta$ be an ample divisor on $X$. When $dim(X) = 1$ and the ground field is $\mathbb C$, the Narasimhan-Seshadri theorem establishes an equivalence between the category of irreducible unitary representations of the fundamental group and the category
of stable bundles with degree $0$. When $dim(X)$ is arbitrary, for bundles with all $c_i = 0$, this equivalence is proved by Mehta and Ramanathan in \cite{mr} using restrictions to curves and stable bundles are realised as irreducible unitary representations of the topological fundamental group.

By an analogue of the Narasimhan-Seshadri theorem we mean establishing an equivalence between the category of irreducible representations of a certain group scheme associated to $X$ and the category
of stable bundles on $X$ when the $c_i$'s do not necessarily vanish. Such an analogue is not known even over fields of characteristic $0$. 
 
The aim of this paper is to define this category, to establish such an equivalence, and to obtain some properties and give a few applications.

To make these considerations more precise and state our first theorem we need a few definitions. A vector bundle is said to be {\em strongly semistable} if all its Frobenius pull-backs are semistable (see Definition \ref{frob}). {\em Over fields of characteristic zero}, the notions of strong semistability and strong stability are the usual notions of Mumford semistability and stability.

We say that a bundle $E$ of degree $0$ is {\em{ locally free graded}} (abbreviated {\em lf-graded})   if it is semistable and has a Jordan-H\"older filtration such that the successive quotients are comprised of stable locally free sheaves (i.e vector bundles). In particular, each stable bundle of degree $0$ is an object in this category. We say that a bundle $E$ of degree $0$ is {\em {strongly lf-graded}} if all its  Frobenius pull-backs $F^l(E)$ are lf-graded (see Definition \ref{strongnewss}).

Our first result is the following (see Theorem \ref{tannaka} and Theorem \ref{ns}): semistability and stability are with respect to the fixed ample divisor $\Theta$ on $X$.
\\
\\
\noindent
{\sc Theorem 1.} ~~ {\it  Let ${\mathcal C}^{\ell f}$ denote the category of strongly lf-graded bundles on $X$ and let $x \in X$. Let $\omega_x:{\mathcal C}^{\ell f} \to Vect_k$ be the fibre functor which sends any bundle $V \in Obj({\mathcal C}^{\ell f})$ to its fibre $V_x$. Then, the pair $({\mathcal C}^{\ell f}, \omega_x)$ forms a neutral Tannaka category. Let $\varpi(X,x)$ be the Grothendieck-Tannaka group scheme associated to 
$({\mathcal C}^{\ell f}, \omega_x)$. Then there exists a universal $\varpi(X,x)$--torsor ${\mathcal E}$ on $X$ with the following properties:
\begin{enumerate}
\item  a representation $\rho: \varpi(X,x) \to GL(n)$ is irreducible if and only if the associated vector bundle ${\mathcal E}(\rho) \simeq {\mathcal E} \times^{\varpi(X,x)} k^n$ is stable. 
\item a bundle $V$ is strongly stable of degree $0$ if and only if $V \simeq {\mathcal E}(\rho)$ with $Im(\rho)_{red}$ {\em irreducible in GL(V)}.
\end{enumerate}}

We call $\varpi(X,x)$ the {\em holonomy group scheme in degree $0$} of the variety $X$ at $x \in X$ (see Definition \ref{vtfg}).

A few words about the method of proof. A well-known property of strongly semistable bundles of degree $0$ is that they are closed under tensor products (cf. Section 3 below). On the other hand, unlike the case of curves, the category of strongly semistable bundles of degee $0$ but with arbitrary higher Chern classes is not an abelian category. 

We show  that the property of lf-gradedness is equivalent to a certain weak restriction property  of semistable bundles to smooth divisors (Definition \ref{wrp}). Using this weak restriction property we prove that the subcategory ${\mathcal C}^{\ell f}$ of strongly lf-graded bundles is an abelian category.

We then go on to prove that the category ${\mathcal C}^{\ell f}$ of strongly lf-graded bundle is closed under tensor products (Proposition \ref{rho}). To prove this we need some results from \cite{rr} and \cite{ch}. Our results are characteristic free and proving the tensor product theorem is a little subtle since {\it it is not known} if restriction of strongly semistable bundles to even general high degree CI curves is strongly semistable and such tensor product theorems usually force representation-theoretic bounds on the characteristic of $k$ (cf. \cite{tens}, \cite{imp}, \cite{serre}, \cite{bapa}).

We define the {\em holonomy group schemes associated to an lf-graded bundle $V$}  as images of the {\em holonomy} representations $\rho$ obtained above and denote it by ${\EuScript H}_{x,\Theta}(V)$.

We have termed our theorem an analogue of the Narasimhan-Seshadri theorem for higher dimensional varieties but unlike the classical theorem, where one ``understands" the category of representations of the fundamental groups, here we have established an equivalence of two categories (namely that of lf-graded bundles and the category of representations of holonomy group schemes) both of which are at present highly intriguing and pose many questions. Nevertheless, we redeem ourselves, by giving a few concrete applications of Theorem 1 which relate the abstract holonomy groups with representations of the fundamental groups of CI curves in $X$. 
\linebreak

\noindent
{\bf Applications}: The first application 
of our Theorem 1 is in Section 5. In this section we work over the ground field $\mathbb C$ and we give a new proof of a result of Balaji and Koll\'ar \cite{balkol}; in fact, we prove the following (effective) generalization of this result. For details, see Theorem \ref{bako} and Theorem \ref{bapako}. 
\\
\\
\noindent
{\sc Theorem 2.} ~~ {\it Let $S(r,c)$ be the set of isomorphism classes of polystable bundles $E$ with $deg(E) = 0$ and such that $rank(E) \leq r$ and $c_2(E) \leq c$. There is a number $\ell = \ell(r,c)$ such that for any $m > \ell$ and {\em any smooth CI curve} $C \in |m \Theta|$, 
\[
{\EuScript H}_{x}(E) = {\EuScript H}_{x}(E|_C)
\]
$\forall E \in S(r,c)$.
 
Further, the group ${\EuScript H}_{x}(E|_C)$ is the Zariski closure of the Narasimhan-Seshadri representation $\pi_1(C,x) \to GL(E_x)$ associated to the bundle $E$}. 

The technique of proof is different from that of \cite{balkol}. The geometry goes into proving Theorem 1 and what remains is an interesting application of some representation theory of finite linear groups. We strengthen the results of \cite{balkol} by making the choices of the CI curves effective and furthermore the result holds for {\em all smooth CI curves} in the linear system. 

We then make some remarks on the Kobayashi-Hitchin correspondence. The Donaldson-Uhlenbeck-Yau theorem shows that a stable bundle is equipped with a canonical Einstein-Hermitian connection (\cite{uy}). Using a result from Biswas (\cite{biswas}, \cite{biswas1}) the  groups ${\EuScript H}_{x, \Theta}(E)$ can be identified with the complexification of the holonomy groups arising from the Einstein-Hermitian connection on $E$ (see Theorem \ref{kohi}). Following Simpson we can view the above picture as a triangle:

\tiny
\beqa
\xymatrix{
\{Stable~bundles\} \ar[dr]_{Theorem 1} \ar[rr]^{D-U-Y}& & 
\{E-H~bundles\} \ar[dl]^{Holonomy~groups} \\
& \{Rep_{\mathbb C}(\varpi(X))\} &
}
\eeqa 
\normalsize
When all the Chern classes of the bundles vanish then in the third vertex we can replace $\varpi(X)$ with the fundamental group $\pi_1(X)$ as in the classical Narasimhan-Seshadri theorem. A natural question which arises from this study is, how close can one get to realizing the Einstein-Hermitian metric by the process of restriction to curves? In Donaldson's approach the restriction to CI curves plays a key role in bounding the Donaldson-Yang-Mills functional.

When $E$ is lf-graded but not stable, the holonomy group can be related to the classical Weil representation of the fundamental groups of curves. We make a few remarks in the context of our theorem (see Theorem \ref{weil-simpson}) relating it to some results of Simpson \cite[Corollary 3.10]{higgs}.

We recall that similar Tannakian constructions have been made earlier by M.V. Nori \cite{nori} where he defines the (true) fundamental group scheme of $X$ by realising it as the Grothendieck-Tannaka group scheme associated to the Tannaka category of {\em essentially finite bundles} (cf. Proposition \ref{truequotient} below).  Nori's group scheme is a natural generalization of the algebraic fundamental group defined by Grothendieck and coincides with it in characteristic $0$. Later, when the ground field is $\mathbb C$, C.Simpson \cite{higgs} studied the Tannaka category of semiharmonic Higgs bundles (hence with all Chern classes vanishing) and realised profound connections between these Tannaka group schemes and non-abelian Hodge theory. The Tannaka category we define is in a sense the natural generalization of these considerations and works  over fields of all characteristics.

The second half of the paper (from Section 6 onwards) is our second application of this construction of the holonomy group scheme. 

Let $G$ be a {\em simple simply connected algebraic group}. We give a construction of $\mu$-stable principal $G$--bundles on smooth projective surfaces $X$; in fact we do more; we construct strongly stable bundles {\em whose holonomy group coincides with $G$} under the assumption that the second Chern class $c_2(E)$ is {\it large} and the characteristic of the field $k$ is larger than the Coxeter number $h_G$ of $G$ ($h_G = {dim(G) \over rank(G)} - 1$) and $k$ is an uncountable algebraically closed field. This solves the existence problem for stable $G$--bundles in positive characteristics and shows the existence of an open subset of stable $G$--bundles whose holonomy group is $G$. Even over fields of characteristic zero this result is new and throws light on the geometry of moduli space of principal bundles on surfaces (cf. \cite{balaji} and \cite{balaji1}).  In particular, this shows that the moduli spaces of principal sheaves constructed in \cite{glss} is {\it non-empty} when $X$ is an algebraic surface. More precisely, we prove the following result (see Theorem \ref{nonemptiness} and Theorem \ref{nonemptinesswithfullhol}):
\\
\\
{\sc Theorem 3.} ~~{\it Let $X$ be a smooth projective surface over an uncountable algebraically closed field $k$ of characteristic $char(k) > h_G$. Let $M_X(G)^s$ denote the moduli space of isomorphism classes of stable principal bundles on $X$ with $c_1 = 0$. Then the set of $k$--valued points $M_X(G)^s(k)$ is non-empty. Furthermore, there exist strongly stable principal $G$--bundles on $X$ whose holonomy group is the whole of $G$}.
\\

The basic strategy is as follows: the first step is to get $SL(2,k)$--bundles on a {\em general plane  curve} with {\it full holonomy} and this involves drawing on the paper \cite{bps}. Using this and some deformation theory, we get such $SL(2,k)$--bundles on the projective plane. Ideas from Donaldson's fundamental paper \cite{donald} play a key role in lifting bundles from curves to surfaces. Finally, by projecting the arbitrary surface $X$ to the projective plane and pulling back stable bundles we construct $SL(2,k)$--bundles $E$ on $X$ which are strongly stable with full holonomy.  In other words, the holonomy representation obtained in the previous theorem surjects onto the structure group of the bundle. The behaviour of the holonomy group scheme under such coverings is a new feature and it plays a key role in what follows as well (see Section 6).

Then we take a principal $SL(2)$ in $G$ whose existence is guaranteed by an extension of Kostant's results on principal 3-dimensional subgroups in semisimple groups to positive characteristics (cf. \cite{balaji}, \cite{balaji1} and \cite{serre}). Using this principal $SL(2)$, we extend the structure group of the $SL(2,k)$--bundle $E$ constructed above to the group $G$. Thus, we get the non-emptiness of the moduli of stable $G$--bundles over surfaces under the assumption that $p > h_G$, the Coxeter number of $G$. We then use deformation techniques as well as the methods developed in Section 6 to prove that there exist strongly stable $G$--bundles whose holonomy group coincides with $G$. Over fields of characteristic zero, these correspond to {\em irreducible} anti-self dual principal $G$--connections (cf. \cite{ahs}).

When $dim(X) = 1$, the existence of strongly stable principal $G$--bundles with full holonomy was shown on general curves in \cite{bps}. The existence of such bundles with full holonomy on an arbitrary smooth projective curve is still not known while in striking contrast we have such bundles on all smooth projective surfaces. The lifting of bundles from  ${\mathbb P}^2$ to arbitrary surfaces coming from techniques developed in this paper help us achieve this.

The problem of non-emptiness becomes quite involved in positive characteristics since representation theory is amenable for study of stable bundles only if the bundles behave well under Frobenius pull-backs (see \cite{rr}). Or else, one may have to impose conditions on the characteristic of the field and even then the stability of the associated bundle is not to be expected. To the best of our knowledge, barring the work of C.Taubes on $4$-manifolds (\cite{taubes}), there is no other strategy of constructing stable $G$--bundles on higher dimensional varieties. However, in the case when $G = GL(n)$, such a non-emptiness result is shown in arbitrary characteristics in \cite{langer2}.

The problem  of construction of stable $G$--bundles remains open for  varieties of dimension $\geq3$. In the context of Hartshorne's conjecture (viz non-existence of stable rank $2$ bundles on ${\bp}^n$, for $n \geq 7$) this therefore could become more delicate. 

{\bf Outline}: The paper is organized as follows. In Section 2, we define and study the category of lf-graded bundles and prove the equivalence of a certain weak restriction property with lf-graded property. In Section 3 we prove that the category of strongly lf-graded bundles is  a neutral Tannaka category. In Section 4 we study  holonomy groups in degree $0$ associated to a strongly lf-graded bundle of degree $0$. We then go on to characterize  irreducible representations of the holonomy group schemes of projective varieties. In section 5 we give some applications of the our main theorem in characteristic zero.

In section 6 we study genuinely ramified morphisms and stability and characterize this property in terms of holonomy group schemes. We close the section by studying the behavior of holonomy group schemes under \'etale morphisms. In sections 7 and 8 we study lifting of stable bundles from curves to surfaces. In Section 9  we construct strongly stable $SL(2,k)$ bundles with full holonomy on general plane curves with the assumption that $k$ is uncountable. In Section 10 we show the existence of principal bundles on the projective plane which are strongly stable. We then complete the proof of the non-emptiness of the moduli space of $G$--bundles on an arbitrary surface. In the last section we make a number of remarks especially in the setting of characteristic zero. In particular, we study the group of connected components of the reductive holonomy group scheme. We conclude the paper with some remarks and questions.

\begin{ack} We thank D.S.Nagaraj, S.M.Bhatwadekar, K.N. Raghavan and M.S.Raghunathan and C.S. Seshadri for some very helpful discussions. We thank Professor P. Deligne for his comments and suggestions. The first author thanks the hospitality of TIFR and the ICTS Programme on Vector Bundles, held in TIFR in March 2008. The second author thanks the hospitality of CMI and IMSc where much of this work was carried out.\end{ack}

\section{Restriction of semistable bundles to divisors}

Let $X$ be a smooth projective variety of $dim(X) = d$ over an algebriacally closed field $k$ with arbitrary characteristic. When $char (k) = p >0$, we define the {\em Frobenius morphism} of $X$ to be $F:X \to X$ such that $F = id_{|X|}$ as a map of topological spaces and on each open set $U \subset X$, $F^{\ast}: {\mathcal O}_X(U)  \to {\mathcal O}_X(U)$ takes $f \to f^p$ for all $f \in {\mathcal O}_X(U)$. In characteristic $0$, we take $F = Id_X$ so that all statements are uniform across characteristics. We recall the following well-known
definition:

\bdefe\label{semistability} Let $\Theta$ be a polarisation on $X$ and define the degree of a torsion-free sheaf $\cf$ to be $deg(\cf) = c_1(\cf) \cdot\Theta^{d-1}$. A torsion-free sheaf $V$ is said to be {\it semistable} (resp {\it stable}) if for every sub-sheaf $W \subset V$, 
\[
 \frac {deg(W)}{rk(W)} \leq \frac {deg(V)}{rk(V)}~~ (resp <).
 \]  \edefe
\noindent

\bdefe\label{frob} Let $E$ be a vector bundle on $X$. We denote by $F^n(E)$ the bundle $(F^n)^*(E)$ obtained by the $n$-fold iterated pull-back of the Frobenius morphism. Define $E$ to be {\em strongly semistable} (resp {\em strongly stable}), if $F^n(E)$ is {\em semistable} (resp {\em stable}) for all $n \geq 1$. \edefe

\brem\label{haha} The above definition can be made for principal $G$--bundles. An easy but important fact that we use repeatedly is the following. Let $E$ be a principal $G$--bundle. Let $F^r_*:G \to G$ be Frobenius homomorphism at the level of groups. Then $F^r(E) \simeq E(F^r_*(G))$, where by  $E(F^r_*(G))$ we mean the associated bundle coming from the homomorphism $F^r_*$. \erem

\subsubsection{\bf lf-graded bundles} We define the basic objects and give its salient properties; these will play a key role in the subsequent sections.

\bdefe\label{newss} We say that a bundle $E$ of degree $0$ is {\em{ locally free graded}} (abbreviated {\em lf-graded})   if it is semistable and has the following property: $E$ has a  filtration  such that the successive quotients are comprised of stable locally free sheaves of degree $0$ (i.e vector bundles). In particular, each stable bundle of degree $0$ is an object in this category.
\edefe

\begin{example} It is well-known that there exist $\mu$--semistable locally free sheaves of rank $2$ with Chern classes $c_1 = 0$ and $c_2 = 1$. Let $X = {\mathbb P}^2$. Then we take generic extensions
\[
0 \to \co \to E \to {\mathcal I}_x \to 0
\]
where ${\mathcal I}_x$ is the ideal sheaf of a point $x \in {\mathbb P}^2$. These extensions are classified by ${\text {Ext}}^1({\mathcal I}_x, \co) = {\mathbb C}$ and a non--trivial extension corresponds to a locally free $\mu$--semistable sheaf which is {\em not} lf-graded.  

\end{example}

\bdefe\label{strongnewss} {\em (The category of strongly lf-graded bundles)} We say that a bundle $E$ of degree $0$ is {\em {strongly lf-graded}} if all its  Frobenius pull-backs $F^l(E)$ are lf-graded. We denote by  ${\mathcal C}^{\ell f}$ the category of strongly lf-graded semistable bundles of degree $0$ on $X$.

\edefe

\blem\label{extensionoflf} Extensions of (strongly) lf-graded  bundles are (strongly) lf-graded.\elem
\pf The proof is immediate from the definitions.\enpf

We recall the effective restriction theorem from 
(\cite[Theorem 5.2]{langer}): Let $\Delta(E) = 2r c_2 - (r-1) c_1^2$ be the discriminant of $E$, $\beta_r = \beta(X,\Theta,r)$ and $m = \Theta^d$ be as in \cite[Section 3]{langer}.
 
\bth\label{bogolang} Let $E$ be a torsion-free sheaf of rank $r \geq 2$ which is stable with respect to $\Theta$. Let $k$ be defined as
\[
k = \lfloor{\frac{r -1}{r} \Delta(E)\cdot {\Theta}^{d-1} +  \frac{1}{mr(r-1)} + \frac{(r-1)\beta_r}{mr}}\rfloor
\]
then $\forall a > k$, and smooth $D \in |a\Theta|$, such that $E|_D$ is torsion-free, then the restriction $E|_D$ is stable on $D$ with respect to $\Theta|_D$.\eeth
\bcor\label{restlfg} Let $E$ be lf-graded  of rank $r$ and degree $0$. Then there exists a constant $k$ such that for any $a > k$, and any smooth divisor $D \in |a\Theta|$ the restriction $E|_{D}$ is lf-graded.\ecor 
\noindent {\sc Proof}. If $E$ is stable, by Theorem \ref{bogolang} such a $k$ exists. Now we induct on rank and assume we have an exact sequence:
\[
0 \to E_1 \to E \to E_2 \to 0
\]
such that $E_2$ is stable and $E_1$ is lf-graded. Then by induction, there is a $k_i$ for $E_i$ then $k = max(k_1,k_2)$ works for $E$. 
\begin{flushright} {\it Q.E.D} \end{flushright}

\brem Observe that the bounds are not invariant under the Frobenius, which is why a strong restriction theorem is difficult in general.\erem

\subsubsection {\bf Weak restriction property} We examine the behavior of semistable bundles under restrictions to complete intersection curves. This will give local criteria for lf-gradedness and enable us to prove that lf-graded bundles form an abelian category.  
\bdefe\label{wrp} Let $E$ be a vector bundle on $X$ and let $p \in X$. The triple $(E,X,p)$ is said to have the {\em {weak restriction property}} (abbreviateed {WRP}) with respect to an ample divisor $\Theta$  if the following recursive property holds: 
\begin{enumerate}
\item $E$ is  semistable of degree $0$ with respect to $\Theta$  and 
\item Given any positive integer $m$, there exists an $a \geq m$, and a smooth divisor $D$, $D \in |a\Theta|$, with $p \in D$, such that the restricted triple $(E|_{D},D,p)$ has WRP on $D$ with respect to the ample divisor $\Theta|_D$ .
\end{enumerate}
\noindent
We say that the bundle $E$ has WRP if the triple $(E,X,p)$ has WRP for every $p \in X$.
\edefe

\brem Note that on curves $(E,C,p)$ has WRP means that the bundle  $E$ is semistable on the curve $C$ of degree $0$ and $p \in C$. \erem

The following proposition is the motivation for making this definition and its proof will need a few lemmas.

\bprop\label{localwrp} Let $X$ be a smooth projective variety of dimension $d \geq 2$ and let $E$ be a bundle such that $(E,X,p)$ has WRP w.r.t $\Theta$. Suppose that we have a quotient of degree $0$
\[
E \to F \to 0
\]
with $F$ stable and torsion-free. Then $F$ is locally free at $p \in X$. \eprop

\blem\label{exactnesswrp} Suppose that we have an exact sequence:
\[
0 \to E_1 \to E \to E_2 \to 0
\]
of bundles of degree $0$ such that $(E,X,p)$ has WRP. Then so do $(E_1,X,p)$ and  $(E_2,X,p)$.\elem

\noindent {\sc Proof}. Immediate from the definition.

\blem\label{localalg} Let $(R, \goth m)$ be a regular local ring of dimension $d$ and suppose that we have an exact sequence of modules:
\[
0 \to N \to M \to T \to 0
\]
where $N$ is a free $R$--module and $M$ is torsion-free and $T$ a torsion $R$--module, and such that codimension of  $Supp(T) \geq 2$. Then, $T = 0$ and $N = M$. \elem

\noindent {\sc Proof}. It is enough to prove this lemma when $M$ is reflexive. To see this, we simply note that for any torsion free module, $M \hra M^{**}$ is an isomorphism outside codimension $\geq 2$. 

The proof is now by induction on $d$. For $d = 2$, $M$ being reflexive is therefore locally free. Now, since both $N$ and $M$ have same ranks, the locus where the inclusion $N \hra M$ is not an isomorphism is a divisor, namely the vanishing of the determinant. Since the codimension of $Supp(T) \geq 2$, it follows that $T = 0$ and $N = M$.

Assume $d \geq 3$. Take a general element $x \in \goth m \subset R$ such that the codimension of $supp(T/xT) \geq 2$. Tensoring the exact sequence by $R/x$, we get 
\[
0 \to N/xN \to M/xM \to T/xT \to 0
\]
is exact since $Tor_1(T, R/xR)$ is a torsion module and $N/xN$ is free. By induction $M/xM$ is free. Since $M$ is a torsion-free $R$--module over a domain $R$, it follows that $x \in \goth m$ is $R$--regular and $M$--regular and hence by
\cite[Lemma 1.3.5]{herzog} we have $proj dim_R(M) = proj dim_{R/xR} (M/xM)$. Since $M/xM$ is free, $proj dim_R(M) = proj dim_{R/xR} (M/xM) = 0$, it implies $M$ is free. Hence again by codimension arguments as above, $T = 0$ and $N = M$.\begin{flushright} {\it Q.E.D} \end{flushright}

\bprop\label{surfacecase} {\em (The case when $X$ is  surface)} Let $X$ be a smooth projective surface. Let $E$ be a bundle such that $(E,X,p)$ has WRP w.r.t $\Theta$. Suppose that we have a quotient of degree $0$
\[
E \to F \to 0
\]
with $F$ stable and torsion-free. Then $F$ is locally free at $p \in X$. \eprop

\noindent {\sc Proof}. Consider the composite map 
\[
E \to F^{**} \to T \to 0
\]
where $T$ is a torsion sheaf supported on the singular locus of $F$. Assume further that, $p \in Supp(T)$. 
Observe that the canonical inclusion $F \hra F^{**}$ is an isomorphism in codimension $1$, i.e the codimension of $Supp(T) \geq 2$. 

Now since $F^{**}$ is reflexive (and therefore locally free since $X$ is a surface) and stable, it follows that any restriction to a smooth divisor $D$ is torsion free. Further, by Theorem \ref{bogolang}, there exists a $k$ such such that for any $a > k$ and every smooth curve $C \in |a\Theta|$ containing $p$, $F^{**}|_C$ is stable. Since $(E,X,p)$ has the weak restriction property, we can choose $C$ as above so that $E|_C$ is semistable. Now consider the map, 
\beqa\label{1}
E|_C \to (F^{**}|_C) \to T|_C \to 0
\eeqa
and let $G = Image(E|_C \to F^{**}|_C)$. 

Moreover, since $E|_C$ is semistable of degree $0$ and $F^{**}|_C$ is stable of degree $0$, $G$ is semistable of degree $0$. Since $G \hra F^{**}|_C$ is map between semistable bundles of same rank and degree $0$ and $F^{**}|_C$ is stable, the map is an isomorphism since $C$ is a curve. Hence, $E|_C$ surjects onto $F^{**}|_C$. Therefore $T|_C = 0$. Hence by Nakayama Lemma, $T \otimes k(p) = 0$ implying $p \notin Sing(F)$, i.e $F$ is locally free at $p$.
\begin{flushright} {\it Q.E.D} \end{flushright}

\noindent
{\sc Proof of Proposition \ref{localwrp}}: The proof is similar to the one given for the case when $X$ is a surface, but there are some important differences. Therefore, at the risk of some repetition, for the sake of clarity we give it in full. 

Consider the composite map 
\[
E \to F^{**} \to T \to 0
\]
where $T$ is a torsion sheaf supported on codimension $\geq 2$. Assume further that, $p \in Supp(T)$. 

Now since $F^{**}$ is reflexive it follows that any restriction to a smooth divisor $D$ is torsion free. Further, by the stability of $F^{**}$, by Theorem \ref{bogolang}, there exists a $k$ such such that for any $a > k$ and every smooth $D \in |a\Theta|$ containing $p$, $F^{**}|_D$ is stable. Since $(E,X,p)$ has the weak restriction property, we can choose $D$ as above so that $(E|_D,D,p)$ has the WRP. Now consider the map, 
\beqa\label{1}
E|_D \to (F^{**}|_D) \to T|_D \to 0
\eeqa
and let $G = Image(E|_D \to F^{**}|_D)$. Then observe that $G$ is torsion-free and hence by induction on dimension and by Proposition \ref{surfacecase}, $G$ is locally free at $p \in D$. 

Further, since $E|_D$ is semistable of degree $0$ and $F^{**}|_D$ is stable of degree $0$, it follows that $G$ is semistable of degree $0$ and the inclusion $G \hra F^{**}|_D$ is an isomorphism in codimension $1$, i.e codimension of $Supp(T|_D) \geq 2$.

Observe that in the exact sequence 
\[
0 \to G \to F^{**}|_D \to T|_D \to 0\
\]
since $G$ is free at $p \in D$, by Lemma \ref{localalg}, $F^{**}|_D$ is free at $p \in D$ and therefore, $(T|_D) \otimes k(p) = 0$. Hence by Nakayama lemma, $T \otimes k(p) = 0$ on $X$.

Since $F^{**}$ is torsion-free, the hypotheses of \cite[Lemma 1.3.5]{herzog} apply to $F^{**}$ in the local ring at $p$. Therefore, since $F^{**}|_D$ is free at $p \in D$, it implies that $F^{**}$ is free at $p \in X$ (see proof of Lemma \ref{localalg} above). This together with the fact that $T \otimes k(p) = 0$ implies that $F$ is free at $p \in X$.
\begin{flushright} {\it Q.E.D} \end{flushright}

\bth\label{weaknori} $E$ is an lf-graded bundle of degree $0$ if and only if $E$ has the weak restriction property.\eeth

\noindent {\sc Proof}. Let $E$ be semistable of degree $0$ which is lf-graded. Then by an immediate application of Lemma \ref{restlfg} and an induction on dimension it follows that $E$ has the weak restriction property. In fact, it has even a much stronger restriction property.

Conversely, let $E$ be semistable of degree $0$ with the weak restriction property. Suppose that we have an exact sequence 
\[
0 \to E_1 \to E \to E_2 \to 0
\]
with stable torsion free quotient $E_2$. By Proposition \ref{localwrp} we see that $E_2$ is locally free since $(E,X,p)$ has WRP for every $p \in X$ .
By Lemma \ref{exactnesswrp}  it follows that $E_i$ have WRP  and are therefore lf-graded by an induction on ranks. Hence by Lemma \ref{extensionoflf}, it follows that $E$ is lf-graded. \begin{flushright} {\it Q.E.D} \end{flushright}

We then have the following key Proposition.

\bprop\label{key} Let $f:V \to W$ be a map of lf-graded bundles on $X$. Then the map $f$ is of constant rank. In particular, $ker(f)$ and $coker(f)$ are lf-graded. \eprop

\noindent {\sc Proof}. The Propostion is obvious when $dim(X) = 1$, since $V$ and $W$ are semistable of degree $0$ on a curve. 

We now induct on dimension. Let $p \in X$ and let 
$r(p)$ be the rank of the map $f_p:V_p \to W_p$. Then, by Lemma \ref{restlfg}, there exists a smooth divisor $D$ containing $p$ such that $V|_D$ and $W|_D$ are lf-graded. Hence by an induction on dimension, $f_p$ has constant rank and since $X$ is irreducible $r(p)$ is constant everywhere. Hence both $ker(f)$ and $coker(f)$ are locally free and by Lemma \ref{exactnesswrp}  have WRP at each point of $X$. The lf-gradedness of kernel and cokernel now follows from Theorem \ref{weaknori}.
\begin{flushright} {\it Q.E.D} \end{flushright}

\section{The Tannaka category of strongly lf-graded bundles}

In the last section we showed that ${\mathcal C}^{\ell f}$ is an abelian category. It is well-known that the tensor product of strongly semistable bundles of degree $0$ is strongly semistable (\cite{rr}). In this section we show that ${\mathcal C}^{\ell f}$ is closed under tensor products.

\subsubsection{\bf Associated bundles and bounding instability degrees}

 Let $G$ be a reductive algebraic group over $k$. Recall that a principal $G$ bundle $E$ over $X$ is said to be {\it
${\mu}$--semistable (resp $\mu$--stable)} if $\forall$ parabolic subgroups $Q$ of $G$,
$\forall$ reduction $\sigma_Q : U \lr E(G/Q)$, for a {\em big open} $U \subset X$ (i.e codimension of $X - U \geq 2$) and $\forall$
ample line bundle $L$ on $G/Q$, the degree of the bundle $\sigma_Q^* E(L)) \geq 0$, where the degree is computed using the polarization $\Theta$.

We remark that it suffices to check the conditions for maximal parabolics $Q$ in which case, since $Pic(G/Q) = {\mathbb Z}$, the choice of $L$ is canonical. 

For a section $\sigma: U \to E(G/Q)$ we will henceforth denote by
$deg(\sigma)$ the number $deg ~\sigma^*(E(L))$. For $GL(n)$ this definition coincides with the one for bundles given in Definition \ref{semistability}.

We recall the following notations from \cite{rr}. 

Let $K$ be a field. Let $G$ be a 
connected reductive group over $K$ which acts on a projective $K$-variety $M$. Let $m$ be a $K$-rational point of $M$ which is not semistable. Let $P(m)$ be its Kempf instability 
parabolic defined over the algebraic closure $\overline{K}$. Recall that by the canonical nature of the Kempf parabolic, if it is defined over the separable closure $K_s$ then it is already defined over $K$. Therefore, $P(m)$ is always defined over a purely inseparable extension of $K$.

We now recall the definition of the $K_s$--scheme $M(P)$ whose ${\overline K}$-valued points are precisely points of the orbit $O_G(m)$ whose associated Kempf parabolic is $P(m)$ (see \cite[Lemma 2.4]{rr}).
Observe that if $M(P)$ has an $L$--valued point for a purely inseparable extension  $L/K$, then $P(m)$ is defined over $L$.

\blem\label{ramram} Let $E$ be a semistable $G$--bundle. Let $\rho:G \to H$ be a representation defined over $k$ with connected component of the center mapping to the center. Let $P \subset H$ be a maximal parabolic subgroup and $x_0 \in C$ denote the {\em generic point}. Then we have: 
\begin{enumerate}
\item If $\sigma$ is a section of $E(H/P)$ such that $\sigma(x_0)$ is semistable with respect to the ample $L$ . Then $deg(\sigma) \geq 0$.

\item  If $\sigma$ is a section of $E(H/P)$ such that $\sigma(x_0)$ is unstable. 
Then there exists a positive integer $N(\sigma)$ such that if $F^N(E)$ is assumed to be semistable for any $N \geq N(\sigma)$ then $deg(F^N(\sigma)) \geq 0$. Hence $deg(\sigma) = \frac{1}{p^N} deg(F^N(\sigma)) \geq 0$.
\end{enumerate}
\elem
\noindent {\sc Proof}. Part (1) is simply \cite[Proposition 3.10(i)]{rr}.
Part (2) follows by combining \cite[Proposition 3.13]{rr} and \cite[Theorem 3.23]{rr}, where the number $N(\sigma)$ is precisely the pure-inseparability degree of the extension $L/K_s$ where the scheme $M(P)$ is defined over $L$.
\begin{flushright} {\it Q.E.D} \end{flushright}

Now we recall the following boundedness result from \cite[Proposition 4.5]{ch}.

\bprop\label{holla} There exists an integer $N$ such that for any $K_s$-rational point $m$ 
of $M$ which is not semistable, the instability flag $P(m)$ is defined over $K_s^{p ^{-N}}$. \eprop

\subsubsection{\bf lf-graded property of associated bundles}

We begin with a few notations. Let $K = k(t)$ the rational function field in one variable. Let $K_s$ be its separable closure and ${\overline K}$ the algebraic closure. Observe that for any smooth projective curve $C$ defined over $k$, the separable closure $k(C)_s$ can de identified with $K_s$ by choosing a finite separable map from $C$ onto ${\bp}^1$. 

\bdefe\label{degreeforGbundles} Let $G$ be an affine algebraic group  and let $E$ be a principal $G$--bundle. We say that $E$ has {\em degree $0$} if for every character $\chi:G \to {\mathbb G}_m$, the associated line bundle $E(\chi)$ has degree $0$ (with respect to the polarization $\Theta$). \edefe 

\brem\label{centertocenter} Let $G$ be an affine algebraic group and let $\rho:G \to GL(V)$ be a finite dimensional $G$--module. Let $E$ be a principal $G$--bundle of degree $0$. Then the associated vector bundle $E(V)$ has degee $0$. Furthermore, we can decompose $V = \oplus V_{i}$ such that $\rho = \oplus \rho_i$ and $\rho_i:G \to GL(V_i)$ maps the center of $G$ to the center of $GL(V_i)$. \erem 

\bth\label{ourholla} Let $\rho:G \to H$ be a representation defined over $k$ with the connected component of center mapping to the center. Then associated to this representation $\rho$ there exists a positive integer $l = l(\rho, K_s)$ with the following property: let $C$ be any smooth projective curve, and let $E$ be a $G$--bundle of degree $0$ on $C$. Then the associated $H$--bundle $E(\rho)$ is {\em semistable of degree $0$} whenever the $l$-th Frobenius pull-back $F^l(E)$ is {\em semistable}. 
\eeth

\noindent {\sc Proof}. Let $E$ be a semistable $G$--bundle on $C$. For $E(\rho)$ to be semistable, we need to check that for any parabolic $P \subset H$ and ample $L$ on $H/P$ and any section $\sigma: C \to E(H/P)$, we must have $deg(\sigma) \geq 0$, where $deg(\sigma) = deg(\sigma^*(E(L)))$.

Now consider $\sigma(x_0) = m$ as a $k(C)$--rational point of $E(H/P)_{k(C)} = M$. Then by Lemma \ref{ramram} there are two possibilities. If $m$ is a semistable point for the action of $G_{k(C)}$ on $M$, then $deg(\sigma) \geq 0$. On the other hand if $m$ is an unstable point in $M$, then the Kempf parabolic $P(m)$ is defined over $K_s^{p^{-N}}$ by Proposition \ref{holla}, where $N$ is dependent only on $\rho$ and $K_s$ and independent of the $k(C)$--rational point $m$ and $C$. By the canonical nature of $P(m)$ this implies by Galois-descent that it is defined over $k(C)^{p^{-N}}$.

Now take $l = N$ and assume that the Frobenius pull-back $F^l(E)$ is semistable. Then the degree of the section $deg F^l(\sigma) \geq 0$. This proves that $deg(\sigma) = \frac{1}{p^l} deg(F^l(\sigma)) \geq 0$.

This implies that if we assume $F^l(E)$ is semistable, then $E(\rho)$ is semistable.
\begin{flushright} {\it Q.E.D} \end{flushright}

\subsubsection{\bf The tensor category}

\blem\label{bholla} Let $\rho:GL(V) \to GL(W)$ be a representation defined over $k$. Then there exists a positive integer $l$ with the following property: for any $GL(V)$--bundle $E$ of degree $0$ on $X$, if the Frobenius pull-back $F^l(E)$ is lf-graded, then associated $GL(W)$--bundle $E(\rho)$ is also lf-graded. \elem

\noindent {\sc Proof}. Choose $l = l(\rho)$ as in Theorem \ref{ourholla}. Observe that by Remark \ref{centertocenter}, we may assume that $\rho$ maps the connected component of the center to the center.
 
We prove this by an induction on dimensions. Let $dim(X) = 2$.
We note that since $F^l(E)$ is lf-graded, by Theorem \ref{bogolang}, there exists $m$  such that for all $a$, with $a > m$, and for any smooth curve $C \in |a\Theta|$, the bundle $F^l(E)|_C$ is semistable of degree $0$. 

Hence by Theorem \ref{ourholla}, $E(\rho)|_C$ is semistable of degree $0$ for all such curves and hence $E(\rho)$ has WRP. Thus by Theorem \ref{weaknori}, $E(\rho)$ is lf-graded on the surface $X$.

Now let $dim(X) = d$ be arbitrary. Then by Lemma \ref{restlfg}, since $F^l(E)$ is assumed to be lf-graded, there exists an $m$ such that for all $a > m$ and smooth divisors $D \in |a\Theta|$, the restriction $F^l(E)|_D$ is lf-graded. Hence by induction on dimension, $E(\rho)|_D$ is lf-graded for all such divisors. That is $E(\rho)$ has WRP. This implies by Theorem \ref{weaknori} that $E(\rho)$ is lf-graded.
\begin{flushright} {\it Q.E.D} \end{flushright}

\bprop\label{rho} Let $E$ be a strongly lf-graded bundle on $X$. Then for any representation $\rho: GL(V) \to GL(W)$  then bundle $E(\rho)$ is strongly lf-graded. \eprop

\noindent {\sc Proof}. Observe that for any $n$, $F^n(E(\rho))$ is also obtained by a representation $F^n(\rho): H \to G$ (by composing $\rho$ with the Frobenius power map, see Remark \ref{haha}). Since $E$ is strongly lf-graded, we now choose $l = l(F^n(\rho))$ as in Theorem \ref{ourholla}. Then since $F^l(E)$ is also lf-graded, by Lemma \ref{bholla} it follows that $F^n(E(\rho))$ is lf-graded. This implies that $E(\rho)$ is strongly lf-graded.   
\begin{flushright} {\it Q.E.D} \end{flushright}

\bth\label{tannaka} Let ${\mathcal C}^{\ell f}$ denote the category of strongly lf-graded bundles of degree $0$ on $X$ as in Definition \ref{strongnewss}. Fix a point $x \in X$. Then the category $({\mathcal C}^{\ell f}, \omega_x)$, where $\omega_x : {\mathcal C}^{\ell f} \to Vect_k$ is the {\em evaluation map} at $x \in X$, is a {\em neutral Tannaka category}. \eeth

\noindent {\sc Proof}. {\it $\cc^{\ell f}$ is closed under tensor products}: This follows by 
Proposition \ref{rho}. {\it $\cc^{\ell f}$ is an abelian category}. This is Proposition \ref{key}. 
\begin{flushright} {\it Q.E.D} \end{flushright}

\bdefe\label{vtfg}{\em (Holonomy group scheme in degree 
 $0$)} We define the holonomy group scheme in degree $0$ of $X$ to be the 
 Grothendieck-Tannaka group scheme $Aut^{\otimes}({\mathcal C}^{\ell f})$ 
 associated to the Tannaka category $({\mathcal C}^{\ell f}, \omega_x)$ and we denote it by ${\varpi} (X,x,\Theta)$.\edefe 

The {\em true fundamantal group scheme} in the sense of Nori is the Tannaka group scheme associated to the category $\cn$ of {\em essentially finite} bundles on $X$ (see \cite{nori} and see \cite[Page 146, 2.34]{delignemilne} for the nomenclature).

\bprop\label{truequotient} The true fundamantal group scheme ${\pi}_1^{true}(X,x)$ is a quotient of $\varpi(X,x.\Theta)$ for any $\Theta$. More precisely, we have a faithfully flat morphism $q: \varpi(X,x,\Theta) \to {\pi}_1^{true}(X,x)$. \eprop

\noindent {\sc Proof}.  We have a canonical functor $F:{\cn} \to {\cc}^{\ell f}$ given in
the obvious way since any {\em essentially finite} bundle in the sense of Nori is strongly lf-graded. In fact, any essentially finite bundle has WRP by the strong definition of semistability given in \cite{nori}. The functor is fully faithful. So we need to only check that any sub-object of an essentially finite bundle within the lf-graded category, is essentially finite. Now any sub or quotient bundle of degree $0$ of an essentially finite bundle is essentially finite by \cite[Proposition 3.7]{nori} and we are done by the criterion for faithful flatness of morphisms of group schemes by \cite[Proposition 2.21]{delignemilne}, or \cite[Proposition 5, Appendix]{nori}] (see also the proof of Theorem \ref{varpicharofgenram} below).
\begin{flushright} {\it Q.E.D} \end{flushright}

\brem It is known that the maximal pro-\'etale quotient of ${\pi}_1^{true}(X,x)$ is the usual \'etale fundamental group ${\pi}_1^{\acute{e}t}(X,x)$ (see \cite{nori} or \cite{delignemilne}). \erem

\brem\label{frobonholonomy} We observe that the Frobenius induces a morphism of holonomy groups: at the level of categories we have $F^*: {\cc^{\ell f}} \to {\cc^{\ell f}}$ given by $V \mapsto F_X^*(V)$ and this gives the map $F_X^*:{\varpi}(X) \to {\varpi}(X)$. \erem


\section{Holonomy group scheme of a strongly lf-graded bundle}

In this section we define and consider representations of the holonomy group scheme associated to strongly lf-graded bundles. We then go on to prove the analogue of the Narasimhan-Seshadri theorem.

\blem Let $V$ be  strongly lf-graded vector bundle. Let $({\mathcal C}(V), \omega_x)$ be the pair where ${\mathcal 
C}(V)$ is the subcategory of ${\mathcal C}^{\ell f}$ generated by vector bundles 
of degree $0$ arising as subquotients of   $\oplus_{a,b} T^{a,b}(V)$, where $T^{a,b}(V):=V^{\otimes 
a}\otimes (V^*)^{\otimes b}$. Then  $({\mathcal C}(V), \omega_x)$ is a 
neutral Tannaka category.
\elem

\noindent {\sc Proof}. By Proposition \ref{key} ${\mathcal C}(V)$ is a {\em full 
subcategory} of $Vect(X)$  which is closed under tensor products and is also an abelian subcategory  of ${\mathcal C}^{\ell f}$. Hence it is a neutral Tannaka category.
\begin{flushright} {\it Q.E.D} \end{flushright}

\bdefe\label{stronglyG} Following Nori, we define a strongly lf-graded principal $G$--bundle as arising from a functor $F: Rep_k(G) \to \cc^{\ell f}$ satisfying the axioms \cite[Page 77]{nori} i.e, $F$ is a strict, exact, faithful tensor functor (cf \cite{higgs}). \edefe

\bdefe\label{thehol} {\em (\em Holonomy group-scheme of a bundle)} Let ${\EuScript H}_{x,\Theta}(V)$ denote the associated Grothendieck-Tannaka group scheme  to the category $({\mathcal C}(V), \omega_x)$. We call ${\EuScript H}_{x,\Theta}(V)$ the holonomy group-scheme of the bundle $V$. 
Let $E$ be a strongly lf-graded principal $G$--bundle on $X$. 
Then, we can define the holonomy group scheme ${\EuScript H}_{x,\Theta}(E)$ associated to $E$ as follows: fix a faithful representation $\rho: G \hookrightarrow GL(V)$. Define ${\EuScript H}_{x,\Theta}(E) := {\EuScript H}_{x,\Theta}(E(V))$.\edefe

\brem Note that by the Tannakian definition of the holonomy group scheme, the structure group of the underlying principal bundle of an lf-graded vector bundle $V$ can always be reduced to ${\EuScript H}_{x,\Theta}(V)$. \erem

\brem\label{frobonholonomy1} As in Remark \ref{frobonholonomy}, we again observe that the Frobenius induces a morphism of holonomy groups of bundles: at the level of categories we have $F^*: {\cc(V)} \to {\cc(F^*(V))}$ given by $W \mapsto F_X^*(W)$ and this gives the map $F_X^*:{\EuScript H}_{x,\Theta}(F^*(V)) \to {\EuScript H}_{x,\Theta}(V)$. Let $F^*:GL(V) \to GL(V)$ be the Frobenius homomorphism of groups. Then as we have seen in Remark \ref{haha}, ${\EuScript H}_{x,\Theta}(F^*(V)) \simeq F^*({\EuScript H}_{x,\Theta}(V))$. Further, being subgroup schemes of $GL(V)$, these are finite type group schemes and hence after sufficiently many Frobenius pull-backs, we get ${\EuScript H}_{x,\Theta}(F^l(V)) \simeq {\EuScript H}_{x,\Theta}(V)_{red}$.\erem

It is now fairly standard (see \cite{delignemilne}) to show that in Definition \ref{thehol} the way we have described ${\EuScript H}_{x,\Theta}(E)$ is independent of the choice 
of the $G$--module $V$. By the choice of the  base point one can non-canonically identify 
${\EuScript H}_{x,\Theta}(E)$ with a  subgroup of $G$. In fact, with some amount of work, one could possibly
characterize ${\EuScript H}_{x,\Theta}(E)$  as the ``smallest subgroup scheme" to 
which the structure group of $E$  can be reduced {\em admissibly}. i.e 
preserving the degree $0$ property.

 \bdefe\label{fullholonomy} We say that a strongly lf-graded 
$G$--bundle $E$ has full holonomy if the  holonomy group scheme 
${\EuScript H}_{x,\Theta}(E)\subset G$ is the whole group $G$ itself. \edefe

\blem\label{strongandfull} A principal $G$--bundle $E$ is strongly 
stable with full holonomy if and only if $E(V)$ is strongly stable 
for every {\em irreducible} $G$--module $V$. In fact if $G\to G'$ is an  
irreducible homomorphism (i.e the image does not lie in any parabolic subgroup) then every strongly stable $G$-bundle with full holonomy 
induces a strongly stable $G'$ bundle. 
\elem

\blem\label{irrep} Let $E$ be a strongly lf-graded bundle on $X$. Then $E$ is strongly stable if and only if the reduced holonomy group ${\EuScript H}_{x,\Theta}(E)_{red} \subset GL(V)$ is an irreducible subgroup (i.e it does not lie in any parabolic subgroup). \elem

\noindent {\sc Proof}. The proofs of both these lemmas are identical. Assume $E(V)$ is strongly stable. By repeated Frobenius pull-backs, it is not hard to see that we may assume that ${\EuScript H}_{x,\Theta}(E)$ is {\em reduced} (by Remark \ref{frobonholonomy1}). Indeed, the holonomy group scheme for a high Frobenius pull-back is the reduced holonomy group scheme of $E(V)$. Suppose that ${\EuScript H}_{x,\Theta}(E) \subset GL(V)$ is not irreducible. Then there exists an ${\EuScript H}_{x,\Theta}(E)$--submodule of $E(V)_x$ which by the definition of ${\EuScript H}_{x,\Theta}(E)$ gives a subbundle of $E(V)$ of degree $0$ contradicting stability of $E(V)$.

Conversely, suppose that ${\EuScript H}_{x,\Theta}(E)_{red} \subset GL(V)$ is {\em irreducible}. Suppose $E(V)$ is strongly lf-graded but not strongly stable.

So assume that there exists a stable torsion-free quotient of $F^*(E(V))$ of degree $0$. This gives rise to a subbundle $F^*(E(V))$ of degree $0$ since the quotient is locally free by Proposition \ref{localwrp}. This gives an ${\EuScript H}_{x,\Theta}(E)$--submodule of $V$ contradicting the irreducibility of ${\EuScript H}_{x,\Theta}(E) \hra GL(V)$.
\begin{flushright} {\it Q.E.D} \end{flushright}

\subsubsection{\bf Stable bundles and irreducible representations of the holonomy group scheme} We consider vector bundles on smooth projective varieties with $deg_{\Theta}(V) = 0$ and arbitrary higher Chern classes. We first have the following lemma:

\blem\label{tautology} There exists a universal ${\varpi} (X,x,\Theta)$--torsor $\mathcal E$ on $X$.\elem

\noindent {\sc Proof}.  Consider the following functor:
\[
Rep_k({\varpi}(X,x,\Theta)) \simeq (\cc^{\ell f}, \omega_x) \hra Vect_X
\]
coming from the natural {\em forget} functor $({\cc}^{\ell f}, \omega_x) \to Vect(X)$. This composite functor satisfies the axioms of a fibre functor given in \cite{nori} implying that we have a universal ${\varpi} (X,x,\Theta)$-torsor $\mathcal E$ on $X$. \begin{flushright} {\it Q.E.D} \end{flushright}

\bth\label{ns}  A vector bundle $V$ of rank $n$ and $deg_{\Theta}(V) = 0$ is strongly lf-graded if and only if 
it arises as ${\mathcal E}(\eta)$, for a rational representation 
\[
\eta:{\varpi} (X,x,\Theta) \to GL(n).
\]
Further, $V \simeq {\mathcal E}(\eta)$ is {\em stable} if and only if $\eta$ is an irreducible representation.
Moreover $V$ is {\em strongly stable} if and only if all Frobenius pull-backs ${(F^n)}^*(\eta)  (= \eta \circ {(F^n_X)}^*$ as in Remark \ref{frobonholonomy}) are {\em irreducible} representations.
\eeth

\noindent {\sc Proof}. It is immediate from the definition of the category $({\mathcal C}^{\ell f}, \omega_x)$ that given any lf-graded bundle $V$, we have an inclusion of categories ${\cc^{\ell f}}(V) \subset {\cc^{\ell f}}$ inducing the representation $\eta$. Hence any such $V$ arises as ${\mathcal E}(\eta)$, where $\mathcal E$ is as in Lemma \ref{tautology}, and conversely.

We need to show that last statement about strong stability.
We first observe that for every $\eta:{\varpi} (X,x,\Theta) \to GL(n)$, the image $Im(\eta)$ can be identified with the holonomy group ${\EuScript H}_{x,\Theta}({\mathcal E}(\eta))$. Note that this is a finite type group scheme (being a subgroup of $GL(n)$). The claim in the theorem now follows from Lemma \ref{irrep}.
\begin{flushright} {\it Q.E.D} \end{flushright}

\section{Relationship with the Narasimhan-Seshadri theorem}

In this section we give the first application of Theorem \ref{ns} when the ground field is ${\mathbb C}$. The first application gives the main theorem of \cite{balkol} as a corollary of Theorem \ref{ns}. In fact, we give an effective version of the theorem proved in \cite{balkol}. 

We then make some remarks on the Kobayashi-Hitchin correspondence. The Donaldson-Uhlenbeck-Yau theorem shows that a stable bundle is equipped with a canonical Einstein-Hermitian connection. The algebraic holonomy groups ${\EuScript H}_{x}(E)$ can be identified with the complexification of the holonomy groups arising from the Einstein-Hermitian connection on $E$. 
We conclude the section with a few remarks on the holonomy groups of semistable (non-polystable) bundles on curves and its relation with a result of C. Simpson which makes our group scheme more amenable for applications.

\subsubsection{\bf Stable bundles in charactersitic zero}
Let $V$ be polystable of degree $0$. In \cite{balkol}, by the process of {\em restriction to high degree 
curves}, the notion of an {\em algebraic holonomy group} of $V$  was defined with some characterizing properties. The aim here is to give a completely new proof of the main theorem of \cite{balkol}; in fact we get an effective version of this result (see Theorem \ref{bapako}) as a consequence of Theorem \ref{ns}. Philosophically, the issues in addressing this result involving geometry are already taken into the proof of Theorem \ref{ns}. So it is natural to expect that the remaining aspects of the proof should essentially arise out of representation theoretic considerations. We show that this is indeed the case.

Recall the following classical theorem due to Jordan (\cite[page 114]{jordan}) (See also \cite{larsen} as well as \cite{weisfeiler} where there is an explicit $J(r) =(r+1)!r^{a~(log~r) + ~b}$ for suitable constants $a,b$.)
 
\bth\label{jordan0}{\em (C.Jordan, 1878)} There is a universal constant $J(r)$ such that if $F \subset GL_r$ is any finite subgroup, then $F$ has an {\em abelian normal subgroup} $N \vartriangleright F$ with $|F/N| \leq J(r)$.\eeth

\noindent {\sc Proof}.(See \cite[Chapter 5, Section 36]{cr} for details). Since $F$ is finite, we can embed $F$ in the subgroup $U_r$ of unitary matrices. It is ingenious then shows that the subgroup $N$ generated by matrices $\{A \mid ||A - I || < 1/2 \} $ is abelian, where $||A|| = trace(A {\overline A^t})$ is the standard norm on $M_r(\mathbb C)$. It is easy to see that $N \vartriangleright F$.

If $A_1, A_2, \ldots A_n$ are the coset representatives of $F/N$, then $||A_i - A_j|| \geq 1/2$, if $i \neq j$. Represent $A$ as a point in ${\mathbb R}^{2r^2}$. Then $||A||^{1/2} = d(0,A)$ in ${\mathbb R}^{2r^2}$.  Since each $A_i$ is in $U_r$, it follows that $d(0,A_i) = \sqrt{r}$. Hence, the $A_i$'s can be seen as points on the surface of  spheres  with radius $\sqrt{r}$. Since $||A_i - A_j|| \geq 1/2$, $d(0, A_i - A_j) \geq 1/{\sqrt{2}}$ so that small spheres of radius  $1/{\sqrt{8}}$ about the $A_i$'s are non-overlapping. Placing these spheres in a shell defined by two spheres of radius $\sqrt{r} \pm  1/{\sqrt{8}}$ and comparing volumes, one can check that the number of such points is bounded in terms of $r$. The bound using this argument for $J(r)$ is $(\sqrt{8r} + 1)^{2r^2} - (\sqrt{8r} - 1)^{2r^2}$ and is due to Schur. \begin{flushright} {\it Q.E.D} \end{flushright}

\bprop\label{jordan1} Let $H \subset GL(r)$ be a reductive subgroup. For any $H$--module $V$ with $dim(V) > J(r)$,  $V$ is not an irreducible module for any finite subgroup $F \subset H$.
\eprop

\noindent {\sc Proof}. Let $V$ be an irreducible representation of $H$ such that $dim(V) > J(r)$. Then by \cite[Corollary, page 62]{serre}, since $dim(V) > |F/N|$, $V$ is not an irreducible representation of $F$ for any finite subgroup $F \subset H$.

\begin{flushright} {\it Q.E.D} \end{flushright}

\blem\label{lovely} Let $H$ be a reductive irreducible subgroup of $GL(V)$ with $dim(V) = r$ and let $H^o$ denote the connected component of identity in $H$.
Assume that $D = {\mathcal D}(H^o)$, the derived group of $H^o$, is non-trivial.  Then the symmetric power $Sym^{J(r)}(V)$ considered as an $H$--module contains an irreducible summand of dimension $\geq J(r)$. 
\elem
\pf Observe that $D \subset H$ is a non-trivial connected semisimple group. Then $V$ as a $D$--module and contains an irreducible submodule $W$ with  $dim(W) > 1$. Then, there is a dominant $\lambda$ on a Borel subgroup of $D$ such that $W \simeq H^0(D/B,L(\lambda))$. Hence it follows that $H^0(D/B, L({J(r)}\lambda)) \subset Sym^{J(r)}(W)$ is a $D$--submodule and $dim(H^0(D/B, L({J(r)}\lambda))) \geq J(r)$. It follows immediately that even as a $H$--module, $Sym^{J(r)}(V)$ contains the required irreducible submodules of dimension at least $J(r)$. \enpf

Let $E$ be a stable bundle and let $H = {\EuScript H}_{x}(E)$ denote the holonomy group of $E$. Similarly, let $H_C = {\EuScript H}_{x}(E|_C)$ whenever $E|_C$ is semistable. Observe that by definition $H_C \subset H$.

\blem\label{lovely1}  Let $E$ be a stable bundle such that $H$ is a finite group, i.e $H^o = (1)$. Then  for any smooth ample curve $C \subset X$, we have  $H_C = H$. \elem

\noindent {\sc Proof}. Consider $Y = E_H$ the $H$-reduction. Then $Y$ is a connected \'etale cover of $X$ with Galois group $H$. Since $C$ is an ample curve in $X$, it follows that $E_H|_C$ is also a connected cover of $C$. On the other hand the holonomy group of this finite bundle is $H_C$. Hence the connectedness of $E_H|_C$ implies that $H= H_C$. 
\begin{flushright} {\it Q.E.D} \end{flushright}

\blem\label{lovely2} Let $E$ be a stable bundle and suppose that $Z(H)^o \neq 1$. Let $C \subset X$ be an ample curve such that $E|_C$ is stable. Then $Z(H_C) = Z(H)$. \elem

\noindent {\sc Proof}. Observe that $H \subset GL(E_x)$ is an irreducible subgroup. Thus, if $Z(H)^o \neq 1$ by Schur's lemma,  $Z(H) = Z(GL(E_x))$, i.e {\em the scalars} are precisely the center of $H$. In this case, $Z(H) \subset H^o$.

On the other hand, characters of $H$ are precisely line bundles of degree $0$ in the Tannaka category ${\mathcal C}(E)$. Now a degree zero
non-torsion line bundle on $X$ restricts to a non-torsion line bundle on $C$
by Lefschetz Theorem, since $Pic (X) \to Pic (C)$ is injective. This gives non-torsion degree $0$ line bundles in the Tannaka category ${\mathcal C}(E|_C)$ i.e a non-torsion character of $H_C$. But if $Z(H_C)^o = 1$, then since $H_C^o$ has no characters, all characters of $H_C$ will be torsion. Hence, $Z(H_C)^o \neq 1$. 

Again since $E|_C$ is also stable, by the reasoning above, it follows that $Z(H_C) = scalars$ and hence $Z(H_C) = Z(H)$.
\begin{flushright} {\it Q.E.D} \end{flushright}

Recall from Theorem \ref{bogolang} the number $k(E)$ which we term the Bogomolov index of $E$:
\[
k(E) = \lfloor{\frac{r -1}{r} \Delta(E)\cdot {\Theta}^{d-1} +  \frac{1}{mr(r-1)} + \frac{(r-1)\beta_r}{mr}}\rfloor
\]

\bprop\label{jordan2} Let $E$ be a stable $GL(V)$--bundle with $H = {\EuScript H}_{x}(E)$. Assume that $D = {\mathcal D}(H^o) \neq 1$. Let $W = Sym^{J(r)}(V)$ and let $E(W) = \oplus A_j$ be a decomposition into stable bundles. Let $k_0 = max_j k(A_j)$ and let $\ell > k_0$. Let $C \in |\ell\Theta|$ be a smooth curve. Then $H_C^o \neq 1$. \eprop

\noindent {\sc Proof}. Observe that by Lemma \ref{lovely}, there is a summand $W_0 \subset Sym^{J(r)}(V)$, with $dim(W_0) > J(r)$ and such that $W_0$ is an irreducible $H$--module. 

By the assumption on the Bogomolov indices $k_0$, $E(W_0)|_C$ is stable. This implies that $W_0$ as an $H_C$--module is irreducible. By Proposition \ref{jordan1} this implies that $H_C^o \neq 1$.
\begin{flushright} {\it Q.E.D} \end{flushright}

\bth\label{bako} Let $X$ be a smooth projective variety and $E$ be a stable $GL(V)$--bundle of rank $r$. Let $V_1 = End(V)$ and $V_2 = S^{J(r)}(V)$. Let $E(V_1) = \oplus_j W_{j}$ and $E(V_2) = \oplus_j A_j$ be the
decomposition of $E(V_i), i =1,2$ with $W_j$'s and $A_j$'s as stable bundles. Let $\ell$ be bigger than the maximum of the Bogomolov indices for all the $W_j$'s and $A_j$'s. Let $C$ be a smooth curve in $ |\ell\Theta| $.
Then ${\EuScript H}_{x}(E) = {\EuScript H}_{x}(E|_C) $.\eeth

\noindent {\sc Proof}. By Lemma \ref{lovely1} we have the theorem when $H^o = 1$ i.e $H$ is finite.

If $H^o \neq 1$ but ${\mathcal D}(H^o) = 1$, i.e $H^o = Z(H)$ then $H/Z(H) = F$ for a finite group $F$. By Lemma \ref{lovely1} and Lemma \ref{lovely2} it follows that $H = H_C$ in this case.

Let $E$ be stable such that ${\mathcal D}(H^o) \neq 1$. Let $C$ be so chosen in $|\ell\Theta|$, with $\ell$ as in the theorem. Then by Proposition \ref{jordan2} if ${\EuScript H}_{x}(E|_C) = H_C$ then $H_C^o \neq 1$.

Consider the adjoint group $H' = H/Z(H)$. Then $Lie(H)$ is an irreducible $H'$--module since it is an irreducible $H^o/Z(H)$--module.

Consider the polystable bundle $E(V_1)$. Note that by the stability of $E$, we have the inclusion $E(Lie(H)) \subset E(V_1)$ as a stable summand and hence $E(Lie(H)) = W_j$ for some $j$. Restricting to $C$  we get $E|_C(Lie(H) = W_j|_C$. 

By the same reasoning, observe that $Lie(H_C) \subset Lie(H)$ is an irreducible submodule for the action of $H_C$. Note that $Lie(H_C) \neq (0)$ since $H_C^o \neq 1$. 

Thus, it follows that  $E|_C(Lie(H_C)) \subset E|_C(Lie(H)) = W_j|_C$ is a non-trivial stable sub-bundle and by the choice of $C$, we know that $W_j|_C$ is also stable. This implies, $E|_C(Lie(H_C)) = E|_C(Lie(H))$, i.e $Lie(H_C) = Lie(H)$. Hence $H^o = H^o_C$. 

Now by considering the associated $H/H^o$--bundle $E_H(H/H^o)$ (a finite bundle) by Lemma \ref{lovely1} we get that $H = H_C$.
\begin{flushright} {\it Q.E.D} \end{flushright}

We can summarise the above results in the following effective result.
Define 
\[
\ell(r,c) =  \lfloor{\frac{t -1}{r} \Delta\cdot {\Theta}^{d-1} +  \frac{1}{mt(t-1)} + \frac{(t-1)\beta_t}{mt}}\rfloor\
\]
with $\Delta = 2tc$ and $t = rank (Sym^{J(r)}(E))$ and $J(r)$ is as in (\ref{jordan0}) above. 

\bth\label{bapako} Let $S(r,c)$ be the set of isomorphism classes of polystable bundles $E$ with $c_1(E) = 0$ and such that $rank(E) \leq r$ and $c_2(E) \leq c$. Let $\ell = \ell(r,c)$. Then for any $m > \ell$ and any smooth curve $C \in |m \Theta|$, $\forall E \in S(r,c)$
\[
{\EuScript H}_{x}(E) = {\EuScript H}_{x}(E|_C)
\] 
Further, the groups ${\EuScript H}_{x}(E|_C)$ are the Zariski closures of the Narasimhan-Seshadri representation $\pi_1(C,x) \to GL(E_x)$ associated to the bundle $E$.   
\eeth

\brem We observe that all the results proven above really do not use the ground field of complex numbers except in interpreting the holonomy group as Zariski closures of the Narasimhan-Seshadri representations. \erem

\brem One can interpret the above theorem as saying that any bounded collection of representations of $\varpi(X,x,\Theta)$ with reductive images can be realised as Zariski closures of the Narasimhan-Seshadri representations of $\pi_1(C,x)$. \erem

\subsubsection{\bf The Kobayashi-Hitchin correspondence} Let $E$ be a polystable bundle on $X$. By \cite{uy}, $E$ admits a unique Einstein-Hermitian connection. Assume for simplicity that $det(E) = {\mathcal O}_X$. Let $\nabla$ be the unique Einstein-Hermitian connection on $E$. Let $SU(r) \subset SL(r)$ be the unitary subgroup and let $E_{SU} \subset E$ be a $C^{\infty}$--reduction of structure group of $E$ to $SU(r)$ so that the corresponding connection on $E_{SU}$ gives $\nabla$ on $E$. This reduction to $SU(r)$ is also unique. Let the $C^{\infty}$--connection on $E_{SU}$ be denoted by $\nabla^{SU}$. Using this connection, by parallel transport at $x \in X$, one can define a Lie subgroup $K \subset SU(r)$, namely the {\em smooth holonomy} group of $E_{SU}$ together with a smooth reduction of structure group of $E_{SU}$ to $K$ (see \cite{Ko}).

We then have the following theorem  (cf. \cite{biswas1}).

\bth\label{kohi}(Biswas) Let $E$ be stable and assume that $det(E) \simeq {\mathcal O}_X$. Let ${\overline K} \subset SL(r)$ be the Zariski closure, in $SL(r)$, of $K \subset SU(r)$. Then we have ${\EuScript H}_{x}(E) \simeq {\overline K}$. In other words, ${\EuScript H}_{x}(E)$ is indeed the complexification of the holonomy group associated to the Einstein-Hermitian connection on $E$.
\eeth

\noindent {\sc Proof}. By  \cite[Proposition 3.5]{biswas}, we get a $C^{\infty}$-reduction of structure group $E_{\overline K} \subset E_{SU}$. Further, $E_{\overline K}$ is preserved by the connection $\nabla^{SU}$ which endows $E_{\overline K}$ with a unique connection inducing $\nabla^{SU}$ on $E_{SU}$. This shows that the $C^{\infty}$-reduction of structure group $E_{\overline K} \subset E$ is in fact a holomorphic reduction of structure group. By the defining property of ${\EuScript H}_{x}(E)$, it follows that ${\EuScript H}_{x}(E) \subset {\overline K}$.

Further, it is shown in \cite{biswas1} that even though $H = {\EuScript H}_{x}(E)$ could be disconnected, it can still be shown that the holomorphic reduction $E_H$ is preserved by the connection $\nabla$. Finally, by the minimality property of $E_{\overline K}$, namely, $E_{\overline K}$ is the minimal reductive reduction preserved by $\nabla$, one concludes that ${\EuScript H}_{x}(E) = {\overline K}$.  
\begin{flushright} {\it Q.E.D} \end{flushright}

\subsubsection{\bf Semistable bundles} We make a few remarks on holonomy group schemes of lf-graded bundles which are not necessarily polystable. 

Let $C$ be a smooth projective curve over a field of characteristic zero and fix a base point $x \in C$. By the considerations above and folklore, it is easy to see that the category of semistable bundles of degree $0$ form a neutral Tannaka category and give a group scheme $\varpi(C,x)$ whose representations give all semistable bundles of degree $0$.

Let $V$ be a semistable bundle of degree $0$ and let ${\EuScript H}_{x}(V) \subset GL(V_x)$ be the associated holonomy  group. As we have remarked, when $V$ is polystable of degree $0$, this group can be realised as the image of the holonomy group scheme or as the Zariski closure of the Narasimhan-Seshadri representation.

By viewing $V$ as a Higgs bundle with trivial Higgs structure, we see that by \cite[Corollary 3.10, page 40]{higgs}, there is a holomorphic (Weil) representation $\rho: \pi_1(C,x) \to GL(V_x)$ which is compatible with the unitary Narasimhan-Seshadri representations giving the canonical polystable bundle $gr(V)$. This is a special Weil representation giving the bundle $V$. The flat connection on the semistable bundle coming as extensions of the Yang-Mills connections on the stable bundles was given by Simpson. So we will call this Weil representation as the Weil-Simpson representation.

 Abstractly, since $V$ is semistable of degree $0$, Weil's theorem applies since the indecomposable summands are of degree $0$. Therefore there are holomorphic representations $\eta$ which realise $V$. Among these, the Weil-Simpson representation $\rho$ is canonical to the extent that it induces the unitary ones on the terms of the associated graded.

Now, viewing $(V,0)$ as a semiharmonic Higgs bundle, and noting that all subbundles   of $T^{a,b}V$ of degree $0$ are simply the sub-local systems, we can  talk of the {\it monodromy subgroup} $M(V,x) \subset GL(V_x)$. This subgroup can be described as a minimal reduction subgroup of $V$. It is easy to see that $M(V,x)$ coincides with ${\EuScript H}_{x}(V)$.

Thus the representation $\rho: \pi_1(C,x) \to GL(V_x)$ factors via $\rho: \pi_1(C,x) \to {\EuScript H}_{x}(V)$ and has the characterizing property that, the quotient $\rho: \pi_1(C,x) \to {\EuScript H}_{x}(V)/ R_u({\EuScript H}_{x}(V))$ is {\it unitary}. 

The total holonomy group scheme $\varpi(C,x)$ can thus be described as the projective limit of the Weil-Simpson type representations of $\pi_1(C,x)$.

\subsubsection{\bf The Simpson connection}\label{simpson}

Recall the following result of Simpson (\cite[Corollary 3.10]{higgs}): there exists an equivalance of categories between holomorphic flat bundles which are extensions of unitary flat bundles and semistable bundles with vanishing Chern classes in degree one and two. Observe that when $c_1 = c_2 = 0$ then any semistable bundle is automatically lf-graded.

We then have the following relationship with ${\EuScript H}_{x,\Theta}(V)$.

\bth\label{weil-simpson} Let $V$ be a semistable bundle on $X$ with $c_1 = c_2 = 0$. Then the holonomy group ${\EuScript H}_{x,\Theta}(V)$  (which is now independent of the polarisation) can be realised as the holonomy group of the canonical Simpson connection on $V$. \eeth

\noindent {\sc Proof}. The Simpson connection is realised as the canonical flat connection obtained as the extension of the unitary connection coming from the summands of the associated polystable bundle. Viewing $(V,0)$ as a Higgs bundle, its monodromy group coming from the canonical flat connection. The rest of the proof is as we have indicated for curves above. \begin{flushright} {\it Q.E.D} \end{flushright}

\brem In this connection and a generalization of Simspon's theorem, we refer the reader to the paper of Biswas-Subramanian \cite{bissub}. \erem

\subsubsection{\bf Some more remarks in Characteristic zero.}
\begin{enumerate}

\item Let ${\mathcal C}^{\bpp}$ be the full subcategory of $\mathcal C^{\ell f}$ consisting of {\em polystable  bundles} of degree $0$. It is easy to see that ${\mathcal C}^{\bpp}$ is also a neutral Tannaka category. Moreover, it is a {\em semisimple} category (see \cite{delignemilne}).

\item Define the 
 {\em pro-reductive quotient}  ${\varpi}(X)_{\bpp}$ of 
 ${\varpi} (X)$ universally as follows: whenever 
 $\rho: {\varpi} (X) \to G$ is {\em surjective} with $G$ reductive,
 the representation $\rho$ factors through a representation 
 ${\rho}':{\varpi}(X)_{\bpp} \to G$. (which is automatically 
 surjective)  

\item The Grothendieck-Tannaka group scheme $Aut^{\otimes}({\mathcal C}^{\bpp})$ is isomorphic to ${\varpi}(X)_{\bpp}$. 

\end{enumerate}

\brem The holonomy group $\varpi(X)$ by its Tannakian definition is an affine group scheme which is realised as an inverse limit of algebraic groups. In particular, the topology that it gets is the {\em projective limit topology}.

With this topology, one can ask the question ``Can one compute $\pi_0(\varpi(X))$?" There are a few examples where we can say something. Let $X$ be a smooth projective variety over characteristic $0$, such that $\pi_1^{\acute{e}t}(X,x) = {1}$. In particular, $X$ has no \'etale covering. Then it is immediate that for every $E \in obj(\cc^{\ell f})$, the holonomy group scheme ${\cm}_{x,\Theta}(E)$ is connected. 

Now to see that the entire group scheme $\varpi(X,x,\Theta)$ is connected, it is enough to check (see \cite[Corollary 2.22]{delignemilne}) that there are no nontrivial epimorphisms to any finite group. If $\phi: \varpi(X,x,\Theta) \to H$ is such a homomorphism, with $H$ finite, then embed $H \hra GL(V)$ and consider the composite $\varpi(X,x,\Theta) \to GL(V)$. This induces an lf-graded bundle with finite and hence disconnected holonomy, contradicting what we have mentioned above. 
\erem

In fact (at least over characteristic $0$), we can say more (cf. \cite{dw}). 

\bprop We have an exact sequence of group schemes:
\beqa\label{conncomps}
1 \to \varpi(X,x,\Theta)_{\bpp}^o \to \varpi(X,x,\Theta)_{\bpp} \to \pi_1^{\acute{e}t}(X,x) \to 1
\eeqa
\eprop

\noindent {\sc Proof}. Observe that the quotient surjection is as in Proposition \ref{truequotient} together with the observation that {\em any bundle which is trivialised by a finite \'etale covering is actually polystable}.

Let $\phi: Y \to X$ be an \'etale Galois covering. This induces a functor $\phi^*:\cc(X)^{\bpp} \to \cc(Y)^{\bpp}$ by taking pull-backs and noting that polystable bundles pull-back to polystable bundles. Thus, by Proposition \ref{varpicharofetale} we have an exact sequence:
\[
1 \to \varpi(Y)_{\bpp} \to \varpi(X)_{\bpp} \to Gal(Y/X) \to 1
\]
By taking inverse limit over Galois coverings of $X$, we get the following exact sequence:
\[
1 \to \underset{\longleftarrow}{\lim}~ \varpi(Y)_{\bpp} \to \varpi(X)_{\bpp} \to \pi_1^{\acute{e}t}(X,x) \to 1
\]
Clearly therefore we have an inclusion $\varpi(X,x,\Theta)_{\bpp}^o \subset \underset{\longleftarrow}{\lim}~\varpi(Y)_{\bpp}$. Hence, to complete the proof we need only show that $\underset{\longleftarrow}{\lim}~\varpi(Y)_{\bpp}$ is connected.

Note that the category of finite dimensional representations of the inverse limit group scheme, $\underset{\longleftarrow}{\lim}~ \varpi(Y)_{\bpp}$, on $k$-vector spaces is the category $\underset{\lr}{\lim}~ \cc(Y)^{\bpp}$. Let $W$ be an object in $\underset{\lr}{\lim}~ \cc(Y)^{\bpp}$. Then we need to show (by \cite[Corollary 2.22]{delignemilne}) that the strictly full subcategory whose objects are isomorphic to subquotients of $n W , n \geq 0$, is {\em not stable} under $\otimes$ (where $n W \simeq \oplus W, (n~copies)$).

Suppose that it is stable under $\otimes$. Then we will show that $W$ is a trivial object in the direct limit category. Since our category is semisimple (because all the categories $\cc(Y)^{\bpp}$ are semisimple) we can get a decompositon of $W$ as:
\[
W = W_1 \oplus \ldots \oplus W_s
\]
where the $W_i$ are simple objects. By assumption, for every $j \geq 1$, the object $W_1^{\otimes j}$ is isomorphic to a subquotient of $n W$ for some $n$. Since we are in characteristic $0$, this implies that $W_1^{\otimes j} = \bigoplus_{i = 1}^{s} m_{ij} W_i$. Let $r > s$ and taking tensor powers  $W_1^{\otimes j}, j = 1, \ldots ,r$, we get an integral dependence relation among the columns of the matrix $(m_{ij})$. 

Now following the argument in \cite[Page 14]{dw}, we get an integral polynomial relation 
\[
P(W_1) = Q(W_1)
\]
Associated to the object $W_1$ we have a chain of bundles $V_i$ on the inverse system of Galois covers $Y_i/X$. Then by the isomorphism above for the polynomials in $W_1$, we get an $l_0$ such that for $l \geq l_0$:
\[
P(V_l) = Q(V_l) ~on~ Y_l
\]
This implies by Weil's theorem that $V_l$'s are finite bundles and hence  trivialised in an \'etale cover of $Y_l$.  Hence, the class of $W_1$ in the category $\underset{\lr}{\lim}~ \cc(Y)^{\bpp}$ is trivial. Similarly all the $W_i$'s are trivial and hence so is $W$. This proves the proposition.
\begin{flushright} {\it Q.E.D} \end{flushright}

\section{Genuinely ramified maps, stability and the holonomy group scheme}

Let $X$ be smooth and projective and $k$ an algebraically closed field of arbitrary characteristic (cf. \cite{paramsub}). In this section we study the behavior of $\varpi(X,x)$ under coverings. This will play a key role in the proof of existence of strongly stable principal bundles. 

Let $E$ be a torsion-free sheaf on $X$. Then one has a {\em unique} filtration called the Harder--Narasimhan filtration,  
$E_{\bullet}:=\{0 = E_0 \subset \cdots \subset
E_l = E\}$ by non-zero sub-sheaves such that each $ gr_i = E_i/E_{i-1}$ is semistable torsion-free and $\mu_i:=\mu(E_i/
E_{i-1})>\mu_{i+1}:= \mu(E_{i+1}/E_i)$. The sub-sheaves $E_i$ are
defined inductively as the inverse image of the maximal sub-sheaf of
maximal slope in $E/E_{i-1}$. The successive quotients $E_i/E_{i-1}$, are termed the
{\em Harder Narasimhan factors} of the sheaf $E$. The sub-sheaf $E_1$ is
called the {\em maximal sub-sheaf} of $E$. This sub-sheaf is {\em semistable} and is denoted by $E_{\rm max}$.
It's slope $\mu(E_1) = \mu(E_{\rm max})$ is called the {\em maximal slope} of
$E$ and denoted by $\mu_{\rm max}(E)$. Notice that one always has: 
 $\mu_{\rm max}(E/E_{\rm max}) < \mu_{\rm max}(E)$.

\blem\label{rank} Let $f: X \to Z$ be a finite separable
morphism of smooth projective varieties.
Then for any semistable vector bundle $W$ on $X$ , $f_*(W)$ is locally free
and we have the inequality
$$\mu_{\rm max}(f_*W)\leq \frac{\mu(W)}{{\rm deg}~f} $$
\elem

\noindent {\sc Proof}. The sheaf $f_*W$ is locally free for any locally free $W$, since $f$ is flat and $Z$ is smooth. The inequality  follows from the fact that
${\text{Hom}}_Z(F,f_*W) \cong {\text{Hom}}_X(f^*F, W)$. Hence semistable bundles of
slope $>\frac{\mu(W)}{{\rm deg}~f}$ have no morphism to $f_*W$.
\begin{flushright} {\it Q.E.D} \end{flushright}

We have the following lemma (cf. \cite[Corollary 1.21]{sommese}).

\blem\label{etalecrit} Let $f:X \to Z$ be a finite separable morphism of smooth projective
varieties. Then $X \to Z$ unramified \'etale if and only if $f_*({\mathcal O}_X)$ is semistable of degree $0$.\elem

\noindent {\sc Proof}. Assume $f_*({\mathcal O}_X)$ is of degree $0$. Its semistability is a trivial consequence of Lemma \ref{rank} above which implies that $\mu_{\rm max} f_*{\mathcal O}_{X}=0$ and the equality is because ${\mathcal O}_Z \subset (f_*{\mathcal O}_{X})_{\rm max}$. Hence $\mu (f_*{\co}_X) = 0$ if and only if
$\mu (f_*{\co}_X) = \mu_{\rm max} f_*{\co}_X$ if and only if $f_*{\co}_X$ is semistable.

Let $R \subset X$ be the ramification locus and $B = f(R) \subset Z$ the branch locus.

Let $\Theta$ be a very ample polarisation on $Z$ and $C$
be a general smooth irreducible complete intersection curve with respect to
$\Theta$. Let $D:=f^{-1}(C)$ be the inverse image of $C$. $C$ (being a general CI curve) will meet the branch locus $B \subset Z$ and hence $D$ will meet the ramification locus $R \subset X$.

By Bertini (applied to the sublinear system coming from the pull-back of the sections of $|\Theta|$ and choosing $C$ to meet $B$ transversally), we see that $D$ is smooth as is the curve $C$. 

Thus, we see that $D \to C$ is \'etale if and only if  
$f:X \to Z$ is \'etale.

Let $I_C$ denote the ideal sheaf of $C$.  Then $I_D=f^*I_C$ is the ideal sheaf of $D$. Now taking the direct image of the exact sequence:
\beqa 
0  \to  I_D  \to  {\co}_X  \to  {\co}_D  \to 0 
\eeqa
we get the exact sequence (since $f$ is finite):
\beqa
0  \to  f_*I_D   \to   f_*{\co}_X    \to  f_*{\co}_D   \to  0 
\eeqa
By the projection formula $f_*I_D = I_C\otimes f_*{\co}_X$, hence we obtain:
\beqa\label{(3)}
0  \to  I_C\otimes f_*{\co}_X \to f_*{\co}_X \to  f_*{\co}_D  \to 0 
\eeqa
Now tensor the exact sequence 
\beqa 
0  \to  I_C  \to  {\co}_Z  \to  {\co}_C  \to 0 
\eeqa
by $f_*{\co}_X$ and comparing with (\ref{(3)}), it follows that $f_*{\co}_D \cong f_*{\co}_X\mid _C$. 

Now observe that $D$ is connected. We argue by induction on $dim(X)$. Assume that $dim(X) = 2$. Now observe that $H^1(X,I_D) = H^1(Z, f_*I_D) = H^1(Z, I_C \otimes f_*{\co}_X) = H^1(Z,(f_*{\co}_X)(-m)) = 0$, by Enriques-Severi Lemma (\cite[Chapter 3, Corollary 7.8]{hartshorne}),  since $m \gg 0$, $C$ being chosen a high degree CI curve. Therefore, $H^0(\co_D) = k$.

Observe that by Mehta-Ramanathan, the semistability of $f_*({\mathcal O}_X)$ will imply that $f_*({\mathcal O}_X)|_C$ is semistable of degree $0$ on $C$. Thus, we see that $f_*{\co}_D$ is semistable of degree $0$ and we are reduced to the case when $dim(X) = dim(Z) = 1$. 

Let $dim(X) = dim(Z) = 1$ and $f_*(\co_X)$ be semistable of degree $0$ on $Z$. By an application of Riemann-Hurwitz formula and \cite[Chapter 4, Ex 2.6]{hartshorne}, we see that $deg(R) = 0$. $R$ being effective, it follows that $f$ is unramified.

The converse is more or less obvious. \begin{flushright} {\it Q.E.D} \end{flushright}

\brem We have in fact the following equivalence:
$f_*{\co}_X$ is semistable $\Longleftrightarrow$ $deg (f_*{\co}_X)  =  0$ $\Longleftrightarrow$ $f$ is \'etale.
\erem

\blem\label{genram}
Let $f:X \to Z$ be a finite separable morphism of smooth projective varieties. Then we have the following:
$(f_*{{\mathcal O}_{X}})_{\rm max}$ forms a sheaf of subalgebras of 
$f_*{\mathcal O}_{X}$ on $Z$ and is also locally free subsheaf of $f_*({\mathcal O}_X)$.
\elem
\noindent {\sc Proof}. Lemma~\ref{rank} shows that $\mu_{\rm max} f_*{\mathcal 
O}_{X}=0$. Again since $f$ is separable, the Harder Narasimhan 
filtration of $f_*{\mathcal O}_{X}$ pulls back to the Harder-Narasimhan filtration and hence we have 
\beqa
\mu_{\rm max} f^*(f_*{\mathcal O}_{X})=0
\eeqa
To prove that $(f_*{{\mathcal O}_{X}})_{\rm max}$ forms a sheaf of subalgebras, we 
need to show that the algebra structure given by the multiplication
\beqa
f_*{{\mathcal O}_{X}}\otimes_{{\mathcal O}_{Z}} f_*{{\mathcal 
O}_{X}}\to f_*{{\mathcal O}_{X}}\eeqa
restricts to a multiplication on $(f_*{{\mathcal O}_{X}})_{\rm max}$. In 
other words we need to show that,
\beqa(f_*{{\mathcal O}_{X}})_{\rm max}\otimes_{{\mathcal O}_{Z}} 
(f_*{{\mathcal O}_{X}})_{\rm max}\to f_*{{\mathcal O}_{X}}\eeqa 
has image contained in $(f_*{{\mathcal O}_{X}})_{\rm max}$. 

Since $(f_*{{\mathcal O}_{X}})_{\rm max}$ is semistable of degree $0$ 
and $\frac{f_*{{\mathcal O}_{X}}}{(f_*{{\mathcal O}_{X}})_{\rm max}}$ 
has maximal slope $<0$, it suffices to prove that $f_*{{\mathcal 
O}_{X}}_{\rm max}\otimes_{{\mathcal O}_{Z}} f_*{{\mathcal O}_{X}}_{\rm 
max}$ is semistable of degree $0$. 
Since each bundle has degree $0$, it suffices to show $(f_*{{\mathcal 
O}_{X}})_{\rm max}\otimes_{{\mathcal O}_{Z}} (f_*{{\mathcal 
O}_{X}})_{\rm max}$ has no sub-sheaf of positive slope. 

Now notice that by the projection formula. 
$$f_*{{\mathcal O}_{X}}\otimes_{{\mathcal O}_{Z}} f_*{{\mathcal 
O}_{X}}\cong f_*f^*(f_*{{\mathcal O}_{X}})$$
Since we have already noted that $f^*(f_*{{\mathcal O}_{X}})$ has 
maximal slope $0$, it has no sub-sheaves of positive slope. 

By Lemma~\ref{rank} we see that $(f_*{\mathcal O}_{X})_{\rm max}$ forms a sheaf of subalgebras of $(f_*{\mathcal 
O}_{X})$. Since $(f_*{\mathcal O}_{X})_{\rm max}$ is torsion-free, it 
is locally free on a big open subset $U$. Taking $Y = {\mathcal 
S}pec(S(f_*{\mathcal O}_{X})_{\rm max}^*)$ and restricting it to $U$, 
we get an \'etale cover $T \to U$ (by Lemma \ref{etalecrit}). 

Since $U \subset Z$ is a big open 
subset, it follows  that $\pi_1^{\acute{e}t}(U) = \pi_1^{\acute{e}t}(Z)$ (by ``purity of branch locus", see\cite[Page 42, Examples 5.2 (h)]{milne}). 

Hence, the 
\'etale cover $T \to U$ extends uniquely to an \'etale cover $g:Y \to Z$.  
It is easy to see that this scheme $Y$ is smooth and $g_*({\mathcal O}_Y) = (f_*{\mathcal O}_{X})_{\rm max}$ implying that $(f_*{\mathcal O}_{X})_{\rm max}$ is locally free.

\begin{flushright} {\it Q.E.D} \end{flushright}

\bdefe\label{genuine} Let $f: X\to Z$ be a finite morphism of smooth
varieties. Then $f$ is said to be {\em genuinely ramified} if $f$ is
separable and does not factor through an \'etale cover of $Z$.
\edefe

\bprop\label{genram1} Let $X$ and $Z$ be smooth projective varieties and 
let $f : X \to Z$ be a finite separable morphism. Then $f$ is genuinely 
ramified if and only if $(f_*{\mathcal O}_{X})_{\rm max}\cong {\mathcal 
O}_Z$.  \eprop

\noindent {\sc Proof}. If $rank ((f_*{\mathcal O}_{X})_{\rm max}) = r > 1$, 
then by Lemma \ref{genram} $f$  factors through $Y$, a non-trivial \'etale cover of rank $r$ and hence it is not genuinely ramified.

Suppose that $f$ is not genuinely ramified and that it 
factors through $g: Y \to Z$ which is unramified. Then $g_*{\mathcal 
O}_Y$ is semistable of degree $0$ (by Lemma \ref{etalecrit}) which is also a subundle of 
$f_*{\mathcal O}_{X}$. Hence $(f_*{\mathcal O}_X)_{\rm max}$
has rank $> 1$. \begin{flushright} {\it Q.E.D} \end{flushright}

\bcor\label{genram2}
Let $X$ be smooth projective variety of $dim(X) = d$. Let $f:X\to {\bp}^d$ 
be a finite separable morphism. Then 
$(f_*{\mathcal O}_{X})_{\rm max}\cong {\mathcal O}_{{\bp}^d}$ i.e $f$ is genuinely ramified.  
\ecor

\bprop\label{stabletostable} Let $f:X\to Z$ be a genuinely
ramified morphism of smooth projective varieties. Then
\begin{enumerate}
\item If $V$ and
$W$ are two semistable bundles on $Z$ of same slope, then
$${\text{Hom}}_Z(V,W)\cong {\text{Hom}}_X(f^*V,f^*W)$$
\item If $V$ is a stable bundle on $Z$, then $f^*V$ is stable on $X$.

\item If $V$ is a semistable bundle
and $W \subset f^*V$ is a subbundle of same slope as $f^*V$,
then $W$ is isomorphic to the pull back of a subbundle of $V$. 
\end{enumerate}
\eprop

\noindent {\sc Proof}. (1) Given two semistable bundles $V$ and $W$ of same slope on
$Z$, we have 
\begin{center}
${\text{Hom}}_X(f^*V,f^*W) \simeq {\text{Hom}}_Z(V, f_*f^*W) \simeq 
{\text{Hom}}_Z(V,W\otimes f_*{\mathcal O}_{X})$. 
\end{center}
Further, since $f$ is genuinely
ramified it follows that $f_*{\mathcal O}_{X}/{\mathcal
O}_{Z}$ has negative maximal slope (see Proposition \ref{genram1}) and we have:
\beqa
{\text{Hom}}_Z(V,W\otimes f_*{\mathcal O}_{X}) \simeq {\text{Hom}}_Z(V,W).
\eeqa
(2) Since the socle (maximal subbundle that is a direct sum stable
bundles (see \cite{mr}) is unique, it follows that the socle of $f^*V$
descends to the socle of $V$ when $f$ is separable. Since $V$ is
stable, this descended bundle has to be $V$ itself. This shows that
the pull back of a stable bundle is polystable under any finite
separable map.  Now the stability of $f^*V$ for genuinely ramified
maps follows from (1) as $f^*(V)$ cannot have endomorphisms.

(3) Let $V$ be a semistable bundle over $C$. Let $W \subset f^*V$ be a
subbundle of same slope. Then the socle $Soc(W)$ of $W$ is contained in
$Soc({f^*V})$ of $f^*V$ and hence a direct summand of
$Soc({f^*V})$. But by uniqueness of the socle, $Soc({f^*V})$ is $f^*(Soc(V))$.
Since stable bundles pull back to stable bundles, $Soc(W)$ coincides with
some factors of $f^*(Soc(V))$ and hence is a pull back. Now the assertion
follows by induction on the rank applied to the bundle $W/{Soc(W)}
\subset {f^*V}/{Soc(W)}$. 

\begin{flushright} {\it Q.E.D} \end{flushright}

\blem\label{induced} Let $f:X \to Z$ be a finite separable morphism of smooth varieties. Let $E$ be an lf-graded bundle on $Z$ with respect to a fixed polarisation on $Z$. Then $f^*(E)$ is lf-graded with respect to the pull-back polarisation. In particular, this induces a homomorphism of group schemes:
\[
f_{*}: \varpi(X,x,f^*(\Theta)) \to \varpi(Z,f(x),\Theta)
\]
\elem

\noindent {\sc Proof}. Observe that if $E$ is stable then $f^*(E)$ is polystable and hence lf-graded. Further, extensions of lf-graded bundles pulls back to extensions of  bundles each of which is lf-graded by induction on ranks. Hence by Lemma \ref{extensionoflf} the lemma follows. 
\begin{flushright} {\it Q.E.D} \end{flushright}   

\bth\label{varpicharofgenram} A finite separable morphism $f:X \to Z$ is {\em genuinely ramified} if and only if induced map of the holonomy group schemes $f_{*}: \varpi(X,x,f^*(\Theta)) \to \varpi(Z,f(x),\Theta)$ is {\em surjective}.\eeth

\noindent {\sc Proof}. By \cite[Proposition 5, Appendix]{nori} (or \cite[Proposition 2.21]{delignemilne}) a map $f:G\to
H$ of affine group schemes is {\em surjective} if and only if the
natural induced functor of Tannaka categories, $f^*: Rep(H)\to Rep(G)$ is {\em fully faithful} and further if
any sequence of $G$-modules $0 \to W' \to f^*V \to W''\to 0$
is obtained by pulling up a sequence $0\to V' \to V \to V'' \to 0$. 

Therefore we see that $f_*:\varpi(Y) \to \varpi(Z)$
is surjective if and only if for every lf-graded bundle $V$ on
$Z$, any lf-graded subbundle $W' \subset f^*V$ is a pull back from an lf-graded bundle on $Z$. 

Assume that $f$ is genuinely ramified. Then the condition follows from Proposition \ref{stabletostable} (3).

Conversely, suppose that $f_*: \varpi(X) \to \varpi(Z)$ is surjective. This implies that the map at the level of categories is fully faithful.

Let the rank of the bundle $(f_*{\mathcal O}_{X})_{\rm max}$ be $r$. This bundle gives a sheaf of algebras as in the proof of Proposition \ref{genram1}. This (by Lemma \ref{etalecrit}) gives rise to an etale cover $g:Y \to Z$ where $g_*{\mathcal O}_Y = (f_*{\mathcal O}_{X})_{\rm max}$. Thus, if $s = H^0({\mathcal O}_X, f^*f_*{\mathcal O}_{X})$, then $s> 1$ if and only
if $r >1$. 

By the full faithfulness, we have an identification of ${\text{Hom}}'s$ and we have therefore:
\[
1= dim ({\text{Hom}}({\mathcal O}_Z, (f_*{\mathcal O}_{X})_{\rm max}) = dim ({\text{Hom}}({\mathcal O}_X, f^*(f_*{\mathcal O}_{X})_{\rm max}))) = s
\] 
because $dim ({\text{Hom}}({\mathcal O}_Z, (f_*{\mathcal O}_{X})) = dim ({\text{Hom}}({\mathcal O}_Z, (f_*{\mathcal O}_{X})_{\rm max}) = 1$. Hence $r = 1$, which   
implies by Proposition \ref{genram1} that the map  $f: X \to Z$ is genuinely ramified. \begin{flushright} {\it Q.E.D} \end{flushright}

We have the following lemma which plays a key role in later applications.

\blem\label{key2} Let $f:X\to Z$ be a genuinely ramified morphism of 
smooth projective varieties. Let $E$ be a principal $G$--bundle on $Z$ 
which is strongly stable with full holonomy. Then, $f^*(E)$ is strongly 
stable $G$--bundle with full holonomy. \elem

\noindent {\sc Proof}. Since $E$ is strongly stable with full holonomy, it 
follows by Lemma \ref{strongandfull} that $E(V)$ is strongly stable for 
every {\em irreducible} $G$ module $V$. We note that this property 
classifies $G$--bundles which are strongly stable with full holonomy. 
By Proposition \ref{stabletostable} (2), $f^*E(V)$ is stable for all irreducible 
$G$-modules, consequently $f^*(E)$ is strongly stable with 
full holonomy. \begin{flushright} {\it Q.E.D} \end{flushright}

\subsubsection{\bf Behavior under \'etale maps} Let $\phi: Y \to X$ be a finite morphism. Then by Lemma \ref{induced}, we have a  homomorphism $\phi_*: \varpi(Y) \to \varpi(X)$.

\blem\label{above1} Let $\phi:Y \to X$ be an \'etale Galois cover. A bundle $W$ on $X$ is lf-graded if and only if $\phi^*(W)$ is lf-graded. \elem

\noindent {\sc Proof}. Let $W \to F \to 0$ be a stable degree $0$ torsion-free quotient. Then we {\em claim} that $F$ is {\em locally free}. Pulling back to $Y$ we get
\[
\phi^*(W) \to \phi^*(F) \to 0
\]
and since $F$ is torsion-free and stable, and $\phi$ is \'etale, $\phi^*(F)$ is torsion-free and polystable (to see this, note that the socle $Soc(\phi^*(F)) \subset \phi^*(F)$ is Galois invariant and hence descends to a subsheaf $F'$ of $F$. Since $F$ is stable, it follows that $F' = F$ and hence $Soc(\phi^*(F)) \simeq \phi^*(F)$.)

Now since $\phi^*(W)$ is assumed to be lf-graded, it has WRP and hence by Proposition \ref{localwrp} we get that $\phi^*(F)$ is locally free. Therefore $\phi_*(\phi^*(F)) \simeq F \otimes \phi_*({\mathcal O}_Y)$ is locally free. Hence, so is $F$ proving the claim. The converse is shown in Lemma \ref{induced}. \begin{flushright} {\it Q.E.D} \end{flushright}

\bcor\label{above2} Let $\phi:Y \to X$ be a finite \'etale morphism. A bundle $W$ on $X$ is lf-graded if and only if $\phi^*(W)$ is lf-graded.\ecor

\noindent {\sc Proof}. Let $Z \to Y \to X$ be the Galois completion and let $\psi:Z \to Y$ and $f:Z \to X$ be the composite. (The existence of such a Galois completion is obvious from classical Galois theory since $X$ is normal. For the general case see \cite[4.4.1.8]{Mr}.) 

Assume that $\phi^*(W)$ is lf-graded. Then,
by Lemma \ref{induced}, $\psi^*(\phi^*(W)) = f^*(W)$ is lf-graded. Hence by Lemma \ref{above1}, $W$ is lf-graded. \begin{flushright} {\it Q.E.D} \end{flushright}

\blem Let $\phi:Y \to X$ be a \'etale Galois cover. Let $W$ be an lf-graded bundle on $Y$. Then $\phi_*(W)$ is lf-graded. \elem

\noindent {\sc Proof}. Consider the diagram (\ref{keydiag}) below, when the map $\phi$ is a Galois covering. Then, it can be regarded as a principal $Gal(Y/X)$--bundle. The pull-back $p:Y \times_X Y \to Y$ always has a canonical section, but since it is a $Gal(Y/X)$--bundle, it implies that $Y \times_X Y \simeq Y \times Gal(Y/X)$.  

Therefore, we see that $\phi^*(\phi_*(W)) \simeq \bigoplus_{g \in Gal(Z/X)} g^*W$ and since $W$ is lf-graded, so is $\phi^*(\phi_*(W))$. This implies by Lemma \ref{above1} that $\phi_*(W)$ is lf-graded. \begin{flushright} {\it Q.E.D} \end{flushright}

\bcor\label{above3} Let $\phi:Y \to X$ be a finite \'etale morphism. Let $W$ be an lf-graded bundle on $Y$. Then $\phi_*(W)$ is lf-graded.\ecor
 
\noindent {\sc Proof}. Again take $Z \to Y \to X$. Since $W$ is lf-graded, so is $\psi^*(W)$ and hence $f_*(\psi^*(W)) \simeq \phi_*(\psi_*(\psi^*(W)))$ is lf-graded. Observe that $W$ is a degree $0$ subbundle of $\psi_*(\psi^*(W)$.

Observe that $deg(\phi_*(W)) =  0$ since $\phi$ is \'etale. Hence, $\phi_*(W)$ is a degree $0$ subbundle of the lf-graded bundle $\phi_*(\psi_*(\psi^*(W)))$.

 This implies that $\phi_*(W)$ has WRP by Lemma \ref{exactnesswrp} and hence is lf-graded. \begin{flushright} {\it Q.E.D} \end{flushright}

\bprop\label{varpicharofetale} Let $\phi: Y \to X$ be an \'etale covering. Then the induced homomorphism $\phi_*: \varpi(Y) \to \varpi(X)$ is a closed immersion. Furthermore, if $\phi$ is also Galois, the we have an exact sequence:
\[
1 \to \varpi(Y) \to \varpi(X) \to Gal(Y/X) \to 1
\]
\eprop

\noindent {\sc Proof}. By \cite[Proposition 2.21 (b)]{delignemilne}, we need to check that if $W$ be an object in $\cc(Y)$. Then, $W$ is isomorphic to a subquotient of an object of the form $\phi^*(W')$ with $W'$ in $\cc(X)$. 

Consider the fibre square

\beqa\label{keydiag}
\xymatrix{
Y \times_X Y \ar[r]_{{p}} \ar[d]_{p} &
Y \ar[d]_{\phi} \\
Y \ar[r]_{\phi} & X 
}
\eeqa

Then, since $\phi:Y \to X$ is \'etale finite, the projection $p: Y \times_X Y \to Y$ is a degree $d$ cover of $Y$, possibly disconnected. Since $\phi$ is flat, it follows that $\phi^*(\phi_*(W)) \simeq p_*(p^*(W))$.
Note that $W$ is a subbundle of $p_*(p^*(W))$  and hence by Corollary \ref{above2} and Corollary \ref{above3} it follows that $W$ is a subbundle of the pull-back of an object in $\cc(X)$, namely $\phi_*(W)$. 

The second half of the proposition follows from the arguments in \cite[Lemma 15]{dw}. \begin{flushright} {\it Q.E.D} \end{flushright}

\bcor Let $\phi:Y \to X$ be a finite \'etale morphism. A bundle $W$ is an lf-graded bundle on $Y$ if and only if $\phi_*(W)$ is lf-graded.\ecor

\noindent {\sc Proof}. One way is simply Corollary \ref{above3}. The converse follows from the proof of Proposition \ref{varpicharofetale} above since $W$ is realised as a degree $0$ subbundle of $\phi^*(\phi_*(W))$ which is lf-graded. \begin{flushright} {\it Q.E.D} \end{flushright}

\brem By Theorem \ref{varpicharofgenram} and Proposition \ref{varpicharofetale}, we get a complete factorization of the induced homomorphism under any finite separable morphism of smooth projective varieties. \erem

\section{Existence of unobstructed stable bundles on a surface}

For this section we assume that the ground field $k$ is an algebraically closed field of  $char(k) \geq 3$ and let $dim(X) = 2$. This section is inspired by some results in Donaldson's paper  (\cite{donald}) where he proves generic smoothness of the moduli space (see also \cite{langer2}).

\bprop\label{newname} Let $M$ be a line bundle on $X$. Then there exists a constant $\alpha(M)$ such that for all $c > \alpha(M)$, there exists a stable $E$ of rank $2$ with $det(E) \simeq \mathcal O_X$ and such $c_2(E) = c$ with the following vanishing property:
\[
h^0(ad(E) \otimes M) = 0.
\]
\eprop

\noindent {\sc Proof}. We note that we can assume, to start with, $det(E) \simeq Q$ with $Q \simeq 2n{\Theta}$, where $\Theta$ is the hyperplane line bundle on $X$. For then, we take $V = E \otimes (-n{\Theta})$. Then $det(V) \simeq {\mathcal O}_X$ and $V$ is also stable and furthermore, $ad(V) \simeq ad(E)$.

\begin{assump}\label{gaga} We now choose $Q$ and $Z$ as follows:

\begin{enumerate}
\item Choose $Q = 2n \Theta$ so that $h^0(Q) > 0$ and $deg(Q) > deg(M)$ so that $h^0(Q^* \otimes M) = 0$.

\item Choose $Z$ so that $H^0(Q \otimes M \otimes I_Z) = 0$. This can for example be made by choosing $Z$ general with ${\ell}(Z) > h^0(Q \otimes M)$. This therefore also implies $h^0(M \otimes I_Z) = 0$ since $Q$ has sections.

\item Choose the length of the cycle ${\ell}(Z)$ as well as the degree of $Q$ (w.r.t $\Theta$) also to be large so that we have stable bundles in the Serre construction (see \cite[Chapter 5]{hl}), i.e
\beqa\label{serre}
0 \rightarrow {\mathcal O}_X \rightarrow E \rightarrow Q \otimes I_Z \rightarrow 0
\eeqa
\end{enumerate}
\end{assump}

We then do the following. Tensor the exact sequence above with the line bundle $M$:
\beqa\label{two}
0 \rightarrow M \rightarrow E \otimes M \rightarrow Q \otimes M \otimes I_Z \rightarrow 0
\eeqa
We need to prove that $H^0(ad(E) \otimes M) = 0$, which follows from the 
following lemma applied to the exact sequence 
\[
 0 \rightarrow ad(E) \rightarrow {\mathcal End}(E) \rightarrow \mathcal O 
\rightarrow 0
\]
where the map ${\mathcal End}(E) \lr \mathcal O$ is the {\it ``Trace"} map.

\blem\label{something} Given $M$, choose $E$ as in Assumption \ref{gaga}. Then  $H^0({\mathcal End}(E) \otimes M) \simeq H^0(M)$. Moreover, any $\phi \in H^0({\mathcal End}(E) \otimes M)$ can be expressed uniquely as $id_E \otimes \psi$ where $\psi \in H^0(M)$. Furthermore, one can identify $\psi = trace(\phi)$. In particular, if $\phi \in H^0(ad(E) \otimes M)$, then $\phi = 0$. \elem

\noindent {\sc Proof}. We first {\it claim} that if $\phi \in H^0(E \otimes M)$ is such that $\phi \circ \theta = 0$, then $\phi = 0$. 

To see this, apply the functor ${\text{Hom}}(-, E \otimes M)$ to the exact sequence (\ref{serre}). Then we get:
\[
0 \rightarrow {\text{Hom}}(Q \otimes I_Z, E \otimes M) \rightarrow {\text{Hom}}(E, E \otimes M) \rightarrow {\text{Hom}}({\mathcal O}_X , E \otimes M)
\]
where the last map is $\phi \rightarrow \phi \circ \theta$. Since $\phi \circ \theta = 0$, it implies that there is a $\gamma \in  {\text{Hom}}(Q \otimes I_Z, E \otimes M)$ which maps to $\phi$. Now any map
\[
\gamma:Q \otimes I_Z \rightarrow E \otimes M
\]
factors via a map $\delta \in {\text{Hom}}(Q , E \otimes M)$ since $E \otimes M$ is locally free and we have a commutative diagram as follows:
\[
\xymatrix{
Q \otimes I_Z \ar[d]_i \ar[r]^{\gamma} & E \otimes M \\
Q  \ar[ru]_{\delta} 
}
\]
where $i: Q \otimes I_Z \hookrightarrow Q$ is the canonical inclusion.
Now by (Assumption \ref{gaga}) we have $h^0(E \otimes M \otimes Q^*) = 0$ (by tensoring the exact sequence (\ref{two}) by $Q^*$). Hence $\delta = 0$, implying $\gamma = 0 = \phi$. This proves our claim.

Let $\phi : E \rightarrow E \otimes M$ and consider the composite 
\[
\phi \circ \theta: {\mathcal O}_X \rightarrow E \otimes M
\]
By (Assumption \ref{gaga}) since $h^0(Q \otimes M \otimes I_Z) = 0$, we have a $\psi: {\mathcal O}_X \rightarrow M$ such that the following diagram commutes:
\beqa
\begin{CD}
  \xymatrix{ 0 \ar[r] & {\mathcal O}_X \ar[r]^{\theta} \ar[d]_{\psi} & E \ar[d]_{\phi} \ar[r] & Q \otimes {I}_Z
    \ar[r] & 0
    \\
    0 \ar[r] & M \ar[r] & E \otimes M \ar[r] & Q \otimes {I}_Z \otimes M \ar[r] & 0\\
}
\end{CD}
\eeqa 
Observe that if we tensor $\psi: {\mathcal O}_X \rightarrow M$ by $E$, we again have a commutative diagram:
\beqa
\begin{CD}
  \xymatrix{ 0 \ar[r] & {\mathcal O}_X \ar[r]^{\theta} \ar[d]_{\psi} & E \ar[d]_{id_E \otimes \psi} \ar[r] & Q \otimes {I}_Z
    \ar[r] & 0
    \\
    0 \ar[r] & M \ar[r] & E \otimes M \ar[r] & Q \otimes {I}_Z \otimes M \ar[r] & 0\\
} 
\end{CD}
\eeqa
where we have $\psi = Trace(id_E \otimes \psi)$. Since $(\theta \otimes id_M) \circ \psi = \phi \circ \theta = (id_E \otimes \psi) \circ \theta$ by the commutativity of the two diagrams, we conclude that:
\[
(\phi - id_E \otimes \psi) \circ \theta = 0
\]
By the claim made above and for suitable choices of $Z$ and $M$, we get:
\[
\phi = id_E \otimes \psi
\]
Hence $\psi = Trace(\phi)$.  This proves the first part of the Lemma. Further, if $\phi \in H^0(ad(E) \otimes M)$, it follows that $\psi = 0$ and hence $\phi = 0$.\begin{flushright} {\it Q.E.D} \end{flushright}

\section{Surjectivity of the restriction map}
In this section the assumptions on the field are as in the previous section namely, $char(k) \geq 3$ with $dim(X) = 2$.

\bth\label{diffsurj} Given a curve $C \subset X$ of genus $g(C) \geq 2$, there exists a constant $\alpha(C)$ such that whenever $c_2 \geq \alpha(C)$, there exists a stable $E$ with $c_2 = c_2(E)$ and such that the restriction map
\[
H^1(adE)) \lr H^1(ad(E|_{C})
\]
is {\it surjective}. In particular, there exists a Zariski open subset $U \subset M_X(SL(2))^{s}$ such that for $E \in U$, the bundle $E|_{C}$ is stable and the restriction map is differentially surjective.\eeth

\noindent {\sc Proof}. Observe that, by Serre duality $H^2(ad(E) \otimes {\mathcal O}_X(-C)) = H^0(ad(E) \otimes {\mathcal O}(C) \otimes K_X)^*$ (Here we use the self-duality of $ad(E)$ since $char(k) > 2$). Therefore, taking $M = {\mathcal O}_X(C) \otimes K_X$  and choosing $E$ as in Proposition \ref{newname}, we get the required vanishing of $H^2$ and the surjectivity. This implies that at the level of infinitesimal deformations we have the {\em surjectivity} of the restriction map of formal schemes:
\[
Def(E) \lr Def(E|_{C}).
\]
Now note that for any bundle on $C$, there are stable bundles in its neighbourhood since the moduli space on curves is irreducible and the stable bundles are dense. Hence we have shown that a general stable bundle on $C$ lifts to a stable bundle on $X$. This proves the theorem.
\begin{flushright} {\it Q.E.D} \end{flushright}

\brem We can see the above deformation argument more transparently 
using stacks as follows. Let ${\mathfrak M}_X(SL(2))$ (resp ${\mathfrak 
M}_C(SL(2))$) be the moduli stack of $SL(2)$--bundle on $X$ (resp $C$). 
The restriction map gives a morphism of stacks from the open substack 
${\mathfrak M}_X(SL(2))^s$ of stable bundles to ${\mathfrak M}_C(SL(2))$. 
The differential of this map at $E \in {\mathfrak M}_X(SL(2))^s(k)$ is 
the map $ H^1(adE) \to H^1(ad(E|_{C}))$ which we have shown to be 
surjective. Hence the image contains a stable bundle on $C$ and the 
differential is surjective at this point too. This proves the required 
surjectivity. \erem 

 \blem\label{mehtaram}{\em (Converse to Mehta-Ramanathan's restriction 
theorem for strongly stable bundles)} Let $P$ be a  $G$--bundle on $X$ 
and let $C \subset X$ be curve belonging to  the 
polarisation $m\Theta$ such that $P|_{C}$ is strongly stable. 
Then $P$ is strongly stable with respect to $\Theta$. \elem

\noindent {\sc Proof}. We  first claim that $P$ is itself a stable
principal $H$--bundle. For, if $Q \subset G$ is a parabolic subgroup
and $\chi$ a dominant character of $Q$, $P_{_Q}$ a $Q$--bundle
obtained from a reduction of structure group to $Q$, note that
\[
deg P_{_Q}(\chi) \cdot m \cdot a = deg P_{_Q}(\chi)_{|_{C}}
\]
where $P_{_Q}(\chi)$ denotes the line bundle associated to the
character $\chi$. Since $a > 0$ it follows by the stability
of $P_{|_{C}}$ that $deg P_{_Q}(\chi)_{|_{C}} > 0$ and hence
$deg P_{_Q}(\chi) > 0$, i.e $P$ is
stable. The Frobenius pull-backs behave similarly. To see this, 
observe that the Frobenius pull back $F^r(P)$ can be realised as 
the associated bundle $P(F^r_*(G))$. Hence if we know strong stability 
on $C$, by observing that taking associated constructions commutes with 
the restriction map, we get the strong stability of $P$ as well. 
\begin{flushright} {\it Q.E.D} \end{flushright}

\section{Holonomy groups of $SL(2)$--bundles on a general plane curve}
Towards constructing bundles on surfaces with full holonomy, we rely on restricting bundles to curves and then lifting back. In this section we construct bundles on plane curves with full-holonomy essentially following \cite{bps}.  
The ground field $k$ is an {\em uncountable algebraically closed field of 
characteristic $p>0$} in 
what follows.

\brem\label{twoelements}Let $G_q := SL(2, {\mathbb F}_q) \subset SL(2,k)$. It is well-known that $G_q$ is generated by the elements 
$\left(\begin{array}{cc}1 & 1 \\0 & 1\end{array}\right)$
and 
$\left(\begin{array}{cc}1 & 0 \\1 & 1\end{array}\right)$

Further, observe that $G_q \subset SL(2,k)$ is an {\em irreducible} subgroup.
\erem

\bprop\label{enriques} Let $C \subset {\mathbb P}^2$ be a general plane curve of genus $\geq 2$. Then there exists a 
strongly stable vector bundle of rank $2$ on $C$ with trivial determinant.
\eprop

\noindent {\sc Proof}.  Choose a nodal plane curve $C_0$ of arithmatic genus 
$g\geq 2$ whose irreducible components are lines in ${\mathbb P}^2$. 
Observe that $g = \frac{(k-1)(k-2)}{2}$, 
where $C_1, \ldots, C_k$ are the irreducible components of $C_0$. 
Note that the \'etale fundamental group of $C_0$ is the profinite 
completion of the free group on $g$ generators. 

Consider the space $S_d$ of degree $d$ curves in ${\mathbb P}^2$. Then 
it is well-known that  $S_d \simeq {\mathbb P}^n$, where $n = {{d+3} 
\choose {3}}$. Hence $S_d$  is {\em irreducible}. 

Assume that $q\, \geq\, 4$.
In Remark \ref{twoelements} we noted that the
subgroup $G_q$ of $G$ is generated by two elements.
 
Since the \'etale fundamental
group of $C_0$ is the profinite completion of the
free group on $g$ generators, with $g \geq 2$, by mapping any two generators to the generators of $G_q$
we have a surjective homomorphism from the
\'etale fundamental group of $C_0$ onto $G_q$.
Hence there is an \'etale Galois covering
of the nodal curve $C_0$ with Galois group $G_q$.

We will prove that
there is a neighborhood $U$ of $C_0$ in $S_d$
such that every curve $C$ in $U$ has a Galois \'etale
covering $D \rightarrow C$
with Galois group $G_q$.

Choose a curve $Spec(R) \rightarrow S_d$ such that the special fibre over $Spec(k) \rightarrow S_d$ is $C_0$ and with general fibre smooth. We thus get a family of curves
$$
C_R \rightarrow \text{Spec} R
$$
which we may assume is a proper separable
morphism of algebraic varieties. Let $K$ denote the
quotient field of $R$ and $C_K$ the generic fiber of $C_R$.
Then we have a {\em specialization} homomorphism
$$
\pi_1^{\acute{e}t}(C_K) \rightarrow \pi_1^{\acute{e}t}(C_0)\, ,
$$
where $\pi_1^{\acute{e}t}$ is the \'etale fundamental group. This
specialization homomorphism is surjective if
the residue field is algebraically closed and if $C_0$
connected (see \cite[9.2, 9.3]{Mr} for a proof).
Now giving a finite Galois
\'etale cover of $C_K$ with Galois group $\Gamma$
is equivalent to giving a surjective group homomorphism
$$
\pi_1^{\acute{e}t}(C_K) \rightarrow \Gamma \rightarrow 1.
$$
If $C_0$ has a finite Galois \'etale cover with Galois group
$\Gamma$, we obtain a surjective group homomorphism
$$
\pi_1^{\acute{e}t}(C_0) \rightarrow \Gamma \rightarrow 1.
$$
Hence if the specialization homomorphism
$\pi_1^{\acute{e}t}(C_K) \rightarrow \pi_1^{\acute{e}t}(C_0)$ is surjective,
then $\Gamma$ is also a quotient of $\pi_1^{\acute{e}t}(C_K)$. In particular, in our case the \'etale cover of the special fibre generizes.

This proves that there is a neighborhood
$U \subset S_d$ of the nodal curve $C_0$
with the property that every curve in the family over
$U \subset S_d$ has a Galois \'etale covering
with Galois group $G_q$.

The covering $D \rightarrow C$ is a Galois \'etale cover with galois group $G_q$ which is therefore a principal $G_q$--bundle. Denote this by $E_D$. Let $E$ be the $SL(2,k)$--bundle obtained from $E_D$ by extension of structure group from the inclusion $G_q \subset SL(2,k)$. Then $E$ is a {\em finite vector bundle} in the language of \cite{nori}. Furthermore, the Frobenius pull-back $F^*(E)$ is obtained by extension of structure group via the composition $G_q \rightarrow SL(2) \rightarrow SL(2)$, where the last map is the Frobenius on the group $SL(2)$. It follows that $F^*(E)$ is also a finite bundle. This implies by \cite[Proposition 3.4]{nori} that $E$ is {\em strongly semistable}.

Moreover, the holonomy of $E$ is precisely $G_q$ which is a reduced irreducible subgroup of $SL(2)$ and hence it follows that $E$ is {\em strongly stable}.
\begin{flushright} {\it Q.E.D} \end{flushright}

\subsubsection{\bf Strongly stable $SL(2,k)$--bundles with full holonomy}

With this Proposition in place, we will now run through some arguments from \cite{bps} which will ensure that there are curves $C \in S_d$ which support strongly stable bundles $E$ such that the holonomy subgroup scheme ${\EuScript H}_{x,\Theta}(E)$ is reduced and in fact coincides with $SL(2,k)$.

\bth\label{bps} {\em (See \cite{bps})} There exist an $SL(2,k)$--bundle $E$  
on a general curve $C \subset {\bp}^2$, which is strongly stable  
with {\em full holonomy group, i.e., ${\EuScript H}_{x,\Theta}(E) = SL(2,k)$.}\eeth

\noindent {\sc Proof}. We outline the proof in the following steps (see \cite[Section 6]{bps})

(1): By Proposition \ref{enriques} we have a general plane curve $C$ and a strongly stable bundle $E$ on it with holonomy group ${\EuScript H}_{x,\Theta}(E) = {\EuScript H}_{x,\Theta}(E)_{red} = G_q$.  Let $H$ be a reductive subgroup of $SL(2,k)$ defined over $\overline{\mathbb F}_p$ . Then for $q \gg 0$, $G_q \nsubseteq H$. Furthermore, none of the conjugates of the $H$ contain $G_q$. Hence
for a fixed $H \subset SL(2,k)$ and its conjugates, we can get an $E$ such that ${\EuScript H}_{x,\Theta}(E) = G_q$ and this $E$ will have no reduction to any conjugate of $H$.

(2): We need an $E$ which is strongly stable with reduced holonomy such that ${\EuScript H}_{x,\Theta}(E)$ is not contained in any reductive $H \subset SL(2,k)$ over $k$. To this end, define the subset:
\begin{center}
$S_{H,n}$ := \{\em {$F^n(E)$ is stable and does
not admit a degree zero reduction to H} \}
\end{center}
\noindent
Then $S_{H,n} \neq \emptyset $. For example, the $E$ gotten above for a large $q$ will lie in $S_{H,n}$. To see this, suppose that some Frobenius power $F^n(E)$ has a degree zero reduction to a $H$. Then, by its defintion there exists an m such that $F^{n+m}(E) \simeq E$. Now if $E_H \subset F^n(E)$ is the degree zero reduction to $H$, then $F^m(E_H) \subset F^{n+m}(E) = E$ will give a reduction of structure group of $E$ which is a contradiction. Indeed, it is a countable intersection of nonempty dense open subsets. 

(3): Define 
\begin{center}
$S := \bigcap_{H,n}$ \{$S_{H,n}$ | $H$ reductive defined over $\overline{\mathbb F_p}$ \} 
\end{center}
and all $n > 0$. By what has been remarked above $S \neq \emptyset$. Let $E \in S$. There is no reason as to why ${\EuScript H}_{x,\Theta}(E)$ is reduced. But we observe that there is always an $r$ such that if $P = F^r(E)$, then the holonomy group of $P$, namely ${\EuScript H}_{x,\Theta}(P)$ is {\em reduced}.

(4): Now observe that by the choice of $E$ and $P$, the subgroup ${\EuScript H}_{x,\Theta}(P) \subset SL(2,k)$ is {\em irreducible} and so is reductive. By the choice of $P$, it follows that ${\EuScript H}_{x,\Theta}(P) = SL(2,k)$ and we are done.  \begin{flushright} {\it Q.E.D} \end{flushright}

\section{Existence of stable principal bundles on surfaces}
In this section the ground field $k$ is an uncountable algebraically closed field of positive characteristic $p \geq 3$. In this section $G$ is a {\em simple simply connected algebraic group}. Recall the notion of a principal three dimensional subgroup of 
 $G$. In characteristic $p$ for $p > h_G$, where 
$h_G$ is the Coxeter number of $G$ (defined as $({dim(G) \over rank(G)} - 1)$, there always exists an irreducible 
subgroup of $G$ (in the sense that it is not contained in any proper parabolic subgroup of $G$) which is the image of the principal homomorphism 
$\rho: SL(2,k) \rightarrow G$ (see \cite{serre}).

\bth\label{principal} 
There exists  strongly stable $SL(2,k)$--bundle $E$  on ${\bp}^2$ such that 
$E|_{C}$ is a strongly stable 
bundle with full holonomy. Furthermore if $p>h(G)$ then the principal 
$G$-bundle $E(\rho)$ is also strongly stable. 
\eeth

\noindent {\sc Proof}. Choose a general plane curve $C$ which supports a strongly stable $SL(2,k)$
bundle with full holonomy which exists by Theorem \ref{bps}. By Theorem \ref{diffsurj} 
we can lift general such bundles on $C$ with full holonomy to a bundle $E$ on ${\bp}^2$. 
Then by Lemma \ref{mehtaram} the bundle $E$ will be strongly stable. Moreover since $E|_C$ has full 
holonomy, this implies that every associated bundle $E(V)|_C$ is strongly stable for all  
irreducible $SL(2,k)$--modules $V$ by Lemma \ref{strongandfull}. By Lemma \ref{mehtaram} it implies that 
$E(V)$ is a strongly stable bundle on ${\bp}^2$. Hence by Lemma \ref{strongandfull} the bundle $E$ has full holonomy.

Since $E|_{C}$ is strongly stable with full holonomy, and since $\rho$ is irreducible 
it implies that $E(G)|_{C}$ is stable. The same argument applies for the Frobenius 
pull-backs proving that $E(G)|_{C}$ is strongly stable. This implies by Lemma \ref{mehtaram} 
that $E(G)$ is strongly stable on ${\bp}^2$.
\begin{flushright} {\it Q.E.D} \end{flushright}

\blem\label{projections} Let $X$ be a smooth projective surface and $\Theta$ be an ample 
line bundle on $X$. Then there exists a genuinely ramified morphism $f:X 
\rightarrow {\mathbb P}^2$ such that $f^*({\mathcal O}(1))$  is some 
power of $\Theta$. \elem

\noindent {\sc Proof}. Embed $X \subset {\mathbb P}^n$ using the  very ample line 
bundle $m\Theta$. Now choose a point in ${\mathbb P}^n$ away from $X$ 
and project. With a little care, one can choose projections 
successively so that the maps are {\it separable}. This can be seen 
as follows: Choose a generic codimension $2$ subspace ${\mathbb P}^{n-2} 
\subset {\mathbb P}^{n}$. This will meet $X$ in a finite set of points 
with multiplicity $1$ each. One can now choose any hyperplane in this 
${\mathbb P}^{n-2}$ which avoids these points. Then projecting from 
this hyperplane we get a map to ${\mathbb P}^{2}$ which is separable.

Consider the resulting map $f:X \rightarrow {\mathbb P}^2$ which is 
chosen to be separable and finite. Since $f$ is obtained by projection it 
has the added property that $f^*{\mathcal O}_{{\mathbb P}^2}(1) = m\Theta$.
\begin{flushright} {\it Q.E.D} \end{flushright}

\bth\label{nonemptiness} There exists an $SL(2,k)$--bundle $P$ on $X$
which is strongly stable with full holonomy group  $SL(2,k)$.
Furthermore, if $SL(2,k) \subset G$ is a principal 
homomorphism (which exists if $p > h_G$), then the associated 
$G$--bundle $P(G)$ is strongly stable with respect to $\Theta$.\eeth

\noindent {\sc Proof}. Let $f:X\to {\mathbb P}^{2}$ be a genuinely ramified 
morphism (obtained by Lemma \ref{projections})  and $E$ on ${\mathbb P}^{2}$ be the vector bundle on ${\mathbb 
P}^{2}$ constructed by Theorem \ref{principal}.
and let $P = f^*(E)$. By Lemma \ref{key2}, $P$ is a strongly stable $SL(2,k)$--bundle with full holonomy. By Lemma \ref{strongandfull}, since $\rho: SL(2,k) \to G$ is a principal homomorphism, its image is an irreducible subgroup of $G$. Hence it follows that $P(G)$ is a strongly stable 
$G$--bundle. This completes the proof of the theorem.
\begin{flushright} {\it Q.E.D} \end{flushright}

\subsubsection{\bf Existence of strongly stable $G$--bundles with full holonomy}
Theorem \ref{nonemptiness} gives the existence of stable $G$ bundles on an arbitrary smooth projective surface. This shows that the moduli space of stable bundles is non-empty. But by this construction, the stable $G$--bundle constructed has its holonomy inside the principal $SL(2,k)$. It is of interest, especially for the geometry of the moduli space of $G$--bundles, to show the existence of stable $G$--bundles whose holonomy group is the whole of $G$ and Theorem \ref{nonemptinesswithfullhol} below proves this result drawing on the techniques developed so far. Even over fields of characteristic $0$, such existence results are not known; these correspond to the {\em irreducible} anti self-dual $G$--connections (cf. \cite[Page 443]{ahs}).

\bprop\label{reducedholonplane} Let $E$ be a strongly stable $G$--bundle on ${\mathbb P}^2$. Then the holonomy group scheme ${\EuScript H}_x(E)$ is {\em reduced and connected}. \eprop

\noindent {\sc Proof}. By Remark \ref{frobonholonomy1}, and the comments following it, we may choose a faithful $G$--module $W$ and work with the strongly semistable vector bundle $V = E(W)$ since ${\EuScript H}_x(E) = {\EuScript H}_x(V)$. Again by Remark \ref{frobonholonomy1}, the Frobenius induces a functor at the level of categories, $F^*: {\cc(V)} \to {\cc(F^*(V))}$ given by $A \mapsto F_X^*(A)$ and this gives the map $F_{{\mathbb P}^{2}}^*:{\EuScript H}_{x}(F^*(V)) \to {\EuScript H}_{x}(V)$. We now claim that: 
\beqa\label{claim}
{\EuScript H}_x(F^*(V)) = {\EuScript H}_x(V). 
\eeqa
To see this, observe firstly that the map $F_{{\mathbb P}^{2}}^*$ is an inclusion ${\EuScript H}_x(F^*(V)) \hookrightarrow {\EuScript H}_x(V)$. To show that this map is surjective, we appeal to \cite[Proposition 2.21(a)]{delignemilne}. The induced functor $F^*: {\cc(V)} \to {\cc(F^*(V))}$ is fully faithful: 
\begin{center}
${\text{Hom}}_{{\mathbb P}^{2}}(F^*U_1,F^*U_2) \simeq {\text{Hom}}_{{\mathbb P}^{2}}(U_1, F_*F^*U_2) \simeq 
{\text{Hom}}_{{\mathbb P}^{2}}(U_1,U_2\otimes F_*{\mathcal O}_{{\mathbb P}^{2}})$. 
\end{center}
On the other hand, $F_*{\mathcal O}_{{\mathbb P}^{2}} = {\mathcal O}_{{\mathbb P}^{2}} \oplus {terms~with~ negative~ degree}$ (cf. for example \cite[Lemma 3.7, page 64]{ein}). 

Since $U_i$ are strongly semistable of degree $0$, it follows easily that ${\text{Hom}}_{{\mathbb P}^{2}}(U_1,U_2\otimes F_*{\mathcal O}_{{\mathbb P}^{2}}) = {\text{Hom}}_{{\mathbb P}^{2}}(U_1,U_2)$. Hence, ${\text{Hom}}_{{\mathbb P}^{2}}(F^*U_1,F^*U_2) \simeq {\text{Hom}}_{{\mathbb P}^{2}}(U_1,U_2)$. 

Let $U \in {\cc(V)}$ and let $A \subset F^*(U)$ be a subobject in ${\cc(F^*(V))}$. So $A$ is a degree $0$ subbundle. In particular, $A$ is $\mu$--semistable. Let $T_{{\mathbb P}^2}$ be the tangent bundle of ${\mathbb P}^2$. We get an ${\mathcal O}_{{\mathbb P}^{2}}$--linear homomorphism $f:T_{{\mathbb P}^2} \to {\mathcal {H}om}(A, F^*(U)/A)$.  Since we are on ${\mathbb P}^2$, $\mu$--semistability is strong semistability and hence $\mu$--semistable bundles are closed under tensor products. Therefore, ${\mathcal{H}om}(A, F^*(U)/A)$ is $\mu$--semistable of degree $0$. Note that $T_{{\mathbb P}^2}$ is semistable with $\mu(T_{{\mathbb P}^2}) > 0$. This implies that $f = 0$. Hence, by standard descent theory for purely inseparable extensions, $A$ descends to a subbundle $A_1 \subset U$ (cf. (\cite[Proposition 1.7]{mehtanori}). This proves the conditions in  \cite[Proposition 2.21(a)]{delignemilne} and hence the claim (\ref{claim}).

As we have already observed (Remark \ref{frobonholonomy1}), by sufficiently large number of Frobenius pull-backs, we see that ${\EuScript H}_x(F^n(V))$ is {\em reduced}. Hence by a repeated application of (\ref{claim}), we get ${\EuScript H}_x(V)$ is {\em reduced}. The simple connectedness of  ${\mathbb P}^2$ now forces that ${\EuScript H}_x(V)$  is also {\em connected}.

\brem It is a fact shown in \cite[Proposition 1.7]{mehtanori} and more generally in \cite{mehtaramanathan} that on smooth projective $X$ such that $\mu(T_X) > 0$, $\mu$--(semi) stability is equivalent to strong (semi) stability. \erem

\bth\label{nonemptinesswithfullhol} Let $X$ be a smooth projective surface and let $char(k) = p$ with $p > h_G$. If  $c_2 \gg 0$ then there exists stable principal $G$--bundles on $X$ whose holonomy group is the whole of $G$. \eeth

\noindent {\sc Proof}. Observe that by Lemma \ref{key2}, if we construct a stable $G$--bundle on ${\mathbb P}^{2}$ with full holonomy, in other words whose holonomy group is the whole of $G$, then by pulling back by a genuinely ramified cover $f:X\to {\mathbb P}^{2}$, we get the required bundle on an arbitrary $X$.

By Theorem \ref{nonemptiness}, the moduli space $M_{{\mathbb P}^{2}}(G, c_2)^s \neq \emptyset$ for large $c_2$ and $p > h_G$. Let $E$ be a stable $G$ bundle on ${\mathbb P}^{2}$. Deformation theory shows that the if $H^2(E(ad G)) = 0$, then the point $E$ is a smooth point in $M_{{\mathbb P}^{2}}(G, c_2)^s$ and further, the component containing $E$ has the expected dimension. By Serre duality, $H^2(E(ad G)) = H^0(E(ad G) \otimes \co(-3))$. Since $E$ is stable, being on ${\mathbb P}^{2}$ it is strongly stable and hence the bundle $E(ad G)$ is semistable of degree $0$ . Therefore $H^0(E(ad G) \otimes \co(-3)) = 0$.

By Proposition \ref{reducedholonplane}, the holonomy groups of strongly stable bundles on ${\mathbb P}^{2}$ are reduced and connected.
Suppose that all the stable $G$--bundles in this smooth component have holonomy inside a smaller reductive group $H \varsubsetneq G$. Now upto conjugacy, there are only countably many reductive subgroups of a semisimple algebraic group $G$. Observe also that conjugate subgroups give isomorphic $G$ bundles in $M_{{\mathbb P}^{2}}(G, c_2)^s$. By a simple dimension count and an application of Baire category theorem, we see that the moduli space $M_{{\mathbb P}^{2}}(G, c_2)^s$ cannot be covered by stable bundles with structure groups in $H \varsubsetneq G$. Hence, there is a stable $G$--bundle on $X$ whose holonomy group is the whole of $G$. \begin{flushright} {\it Q.E.D} \end{flushright}

\brem Existence of a smooth point on the moduli space of $G$ bundles on an arbitrary surface is still unknown. This is true on ${\mathbb P}^{2}$ which is what we use. \erem  

\brem In the case when $G = GL(n)$, such non-emptiness results are shown over fields of arbitrary characteristics in \cite{langer2}. \erem

\section{Miscellaneous remarks}

\brem All the Tannakian constructions in Sections 1-4 go through for big open subsets (complement codim $\geq 2$) $U \subset X$. Furthermore, $\varpi(U,x,\Theta)$ is well-defined and there is a natural map $({\cc}^{\ell f}_X, x) \hra ({\cc}^{\ell f}_U,x)$ inducing a homomorphism $\varpi(U,x) \to \varpi(X,x)$. Since morphisms of vector bundles extend across codimension $\geq 2$, it follows by the surjectivity criterion on Tannaka categories that this homomorphism is surjective. \erem

\brem Notice that $\varpi({\bp}^1)$ is {\em trivial}. On the other hand, $\varpi({\bp}^1 \times {\bp}^1)$ is non-trivial by Theorem \ref{nonemptiness}.  \erem

\brem $\varpi({\bp}^2)$ has no characters since degree $0$ line bundles on $\bp^2$ are trivial. If $Z \to {\bp}^2$ is the blow-up at a point, then $Z$ has the rank of the Neron-Severi group $NS(Z)$ is $2$ implying existence of degree $0$ line bundles. Thus $\varpi(Z)$ has non-trivial characters. This implies that even as abstract groups $\varpi(Z) \neq \varpi({\bp}^2)$. This shows that $\varpi$ is not a birational invariant. \erem 

\brem We have remarked that ${\varpi} (X)$ depends on the choice of the polarisation. It will be interesting to see how the group
scheme changes with the polarisation and what happens with the wall phenomenon.
\erem


\brem Throughout this paper, we assume that degree of all our bundles is $0$. In the light of \cite{bispar}, one should be able to define the holonomy group scheme for strongly lf-graded bundles with arbitrary $\mu$. \erem

\brem After an earlier version of this paper was posted in the archives, we received a manuscript from Adrian Langer \cite{langersfunda}, where he studies the Tannaka category of strongly semistable bundles with all Chern classes zero (analogous to the semiharmonic Higgs bundles). Langer informed us that  his construction is the ``zeroth graded piece" in the structure discussed in this paper.\erem

\subsubsection{\bf Some remarks on the graded Tannaka structure.} These remarks follow from a suggestion of Nori. Let ${\cc}^{\ell f}_o$, be the subcategory of $\cc^{\ell f}$ consisting of bundles with $c_1 = 0$. This allows us to define the following {\em secondary slope}:
\[
{\mu}_2(E) := \frac{c_2(E) \cdot \Theta^{d-2}}{rank(E)}
\]
\blem Let $V$ and $W$ be bundles with $c_1 = 0$. Then
\[
{\mu}_2(V \otimes W) = {\mu}_2(V) + {\mu}_2(W)
\]
\elem

\noindent {\sc Proof}. We see this easily as follows: the Chern character is given by
\[
ch(E) = r + c_1(E) + {1\over2}(c_1^2 - 2 c_2) + \ldots  
\] 
It satisfies $ch(V \otimes W) = ch(V) ch(W)$. Since $c_1 = 0$, we have the equation
\[
{ch(V \otimes W) \over rs} = {ch(V) \over r} {ch(W) \over s}
\]
implying
\[
1 - {c_2(V \otimes W) \over rs} + \ldots = (1 - {c_2(V) \over r} + \ldots) (1 - {c_2(W) \over s} + \ldots).
\]
Hence, comparing terms degreewise, the formula for $\mu_2$ follows.

\end{document}